%% file: main.tex
\documentclass[11pt,a4paper]{article}

\usepackage[american]{babel}
\usepackage[utf8]{inputenc}
\usepackage[T1]{fontenc}

\usepackage[a4paper,top=2.5cm,bottom=2.5cm,left=2.5cm,right=2.5cm,marginparwidth=2cm]{geometry}

\usepackage{amsmath}
\usepackage{graphicx}
\usepackage[colorinlistoftodos]{todonotes}
\usepackage[colorlinks=true, allcolors=blue]{hyperref}
\usepackage{lipsum} 
\usepackage{floatrow} 
\usepackage[hypcap=false]
{caption} 
\usepackage{subcaption}

\usepackage{listings}
\usepackage{url}
\usepackage{courier}
\usepackage{wrapfig}
\usepackage{color,soul}
\definecolor{codegreen}{rgb}{0,0.6,0}
\usepackage{chngpage}
\usepackage{xcolor}
\usepackage{bm}
\usepackage{amsmath,amssymb}
\usepackage{verbatim}
\newlength{\mycolwidth}
\usepackage{array}
\newcolumntype{Z}{>{$}p{\mycolwidth}<{$}}
\usepackage{algorithmicx}
\usepackage{algorithm} 
\usepackage{refcount}
\usepackage{algpseudocode} 
\usepackage{placeins} 
\MakeRobust{\Call}

\usepackage{amsfonts}
\usepackage{graphicx}
\usepackage{epstopdf}
\usepackage{amsmath} 
\usepackage{booktabs}
\usepackage{comment}

\usepackage{stmaryrd}
\renewcommand{\d}{\mathrm{d}}

\DeclareMathAccent{\svec}{\mathord}{letters}{126}
\newcommand\stvec[1]{\mathbf #1}				
				
\newcommand\ssvec[1]{\svec{\stvec{#1}}}	
	
\newcommand\cssvec[1]{\svec{\tilde{\stvec{#1}}}} 
 
\usepackage{authblk}
\usepackage[labelformat=simple]{subcaption}
\usepackage{cite}

\usepackage{lineno}

\definecolor{fortrankeyword}{RGB}{0,0,150}
\definecolor{fortrancomment}{RGB}{120,120,120}
\definecolor{fortranstring}{RGB}{0,128,0}
\definecolor{fortranpragma}{RGB}{170,0,0}
\definecolor{backgray}{RGB}{248,248,248}

\lstdefinestyle{fortranACC}{
  language=Fortran,
  basicstyle=\ttfamily\tiny,
  backgroundcolor=\color{backgray},
  keywordstyle=\color{fortrankeyword}\bfseries,
  commentstyle=\color{fortrancomment}\itshape,
  stringstyle=\color{fortranstring},
  frame=single,
  rulecolor=\color{black!30},
  numbers=left,
  numberstyle=\tiny\color{gray},
  stepnumber=1,
  numbersep=8pt,
  breaklines=true,
  breakatwhitespace=true,
  captionpos=b,
  tabsize=4,
  keepspaces=true,
  showstringspaces=false,
  morecomment=[l]{!\$},                   
  moredelim=[l][\color{fortranpragma}\bfseries]{!\$} 
}

\usepackage{float}
\floatstyle{plain}
\newfloat{listing}{tbp}{lop}
\floatname{listing}{Listing}

\title{HORSES3D-GPU: A high-order discontinuous Galerkin solver for multi-GPU systems}
\author[1,3]{Gerasimos Ntoukas\thanks{Corresponding Author\newline \hspace*{1.55em} E-mail address: gerasimos.ntoukas@upm.es (G. Ntoukas)}}
\author[1,2]{Gonzalo Rubio\thanks{Corresponding Author\newline \hspace*{1.55em} E-mail address: g.rubio@upm.es (G. Rubio)}}
\author[1]{Abbas Ballout}
\author[1]{Stefano Colombo}
\author[1]{David Huergo}
\author[1]{Eduardo Jané}
\author[1]{Albert Jiménez-Ramos}
\author[1]{Hatem Kessasra}
\author[1]{Himpu Marbona}
\author[1]{Oscar Mariño}
\author[1]{Rodrigo Salado}
\author[1]{Pol Solé-Miró}
\author[1]{Andrés M. Rueda-Ramírez}
\author[1,2]{Miguel Chávez-Módena}
\author[1,2]{Eusebio Valero}
\author[1,2]{Esteban Ferrer}
\affil[1]{ETSIAE-UPM - School of Aeronautics, Universidad Politécnica de Madrid, Plaza Cardenal Cisneros 3, E-28040 Madrid, Spain}
\affil[2]{Center for Computational Simulation, Universidad Politécnica de Madrid, Campus de Montegancedo, Boadilla del Monte, 28660 Madrid, Spain}
\affil[3]{Simerics Inc, Bellevue, WA, United States}

\date{}

\begin{document}

\maketitle
\thispagestyle{empty}

\begin{abstract}


We present the GPU acceleration and large-scale performance assessment of \texttt{HORSES3D}, an open-source high-order discontinuous Galerkin solver for computational fluid dynamics. The solver is ported to NVIDIA GPU architectures using OpenACC directives, preserving the original Fortran code structure while enabling GPU-resident execution of the main computational kernels. The implementation exploits the element-local structure of discontinuous Galerkin spectral element methods by mapping element-level loops to GPU gangs and nodal operations to vector-level parallelism.

The GPU version is verified using the method of manufactured solutions and validated on canonical turbulent-flow benchmarks. Its performance is assessed on the MareNostrum~5 accelerated partition using NVIDIA H100 GPUs. Taylor--Green vortex benchmarks show that solver efficiency improves with polynomial order and that near-ideal strong and weak scaling is obtained when the workload exceeds approximately $16,000$ to $20,000$ elements per GPU.

The solver is further evaluated on the High-Lift Common Research Model wing--body configuration, which involves a complex geometry, realistic boundary conditions, and unstructured meshes with (up to 20.8 million) hexahedral elements. Simulations with polynomial orders up to $P=7$ reach approximately $10.7\times10^9$ degrees of freedom and scale efficiently to 2048 GPUs. The results demonstrate that HORSES3D preserves its performance characteristics for industrially relevant configurations and can exploit modern GPU-based supercomputers for billion-degree-of-freedom high-order CFD simulations.

\end{abstract}
\mbox{}
\vfill
\newpage
\input{Methodology}

\bibliographystyle{vancouver}
\bibliography{biblio}
\end{document}

%% file: Methodology.tex


{\bf PROGRAM SUMMARY}

\begin{small}
\noindent
{\em Program Title:} {\fontfamily{qcr}\selectfont HORSES3D-GPU} \\
{\em CPC Library link to program files:} (to be added by Technical Editor) \\
{\em Developer's repository link:}
\url{https://github.com/horses-framework/HORSES3D-gpu}
(GPU version); \url{https://github.com/horses-framework/HORSES3D}
(main framework). \\
{\em Code Ocean capsule:} (to be added by Technical Editor) \\
{\em Licensing provisions:} MIT License \\
{\em Programming language:} Fortran 2008 with OpenACC directives \\
{\em External routines/libraries:} MPI, HDF5 and METIS; the NVIDIA HPC SDK
is used to compile the OpenACC GPU implementation. MKL and PETSc are
optional for CPU and implicit-solver configurations. \\

{\em Nature of problem:}
{\fontfamily{qcr}\selectfont HORSES3D-GPU} is a high-order discontinuous
Galerkin framework for computational fluid dynamics on multi-GPU systems.
The program solves the compressible Navier--Stokes equations for laminar,
transitional and turbulent flows, including configurations with shocks,
complex geometries and wall-bounded turbulence. It supports direct and
large-eddy simulations, several subgrid-scale and wall models, and a broad
selection of numerical fluxes. The wider {\fontfamily{qcr}\selectfont
HORSES3D} framework also includes incompressible and multiphase flows,
immersed-boundary and actuator-line methods, particle dynamics and
aeroacoustics. \\

{\em Solution method:}
The spatial discretisation is based on the nodal discontinuous Galerkin
spectral element method (DGSEM) on curvilinear hexahedral meshes, using
Gauss--Legendre or Gauss--Lobatto--Legendre collocation points. Standard
weak-form and split-form flux-differencing formulations are available,
together with several approximate Riemann solvers. Viscous terms are
discretised using the Bassi--Rebay~1 scheme. The solution is advanced using
explicit low-storage Runge--Kutta methods. The framework supports
arbitrary polynomial orders and $h/p$ discretisations. \\

{\em Additional comments including restrictions and unusual features:}
The GPU implementation uses OpenACC to preserve a single Fortran source
code for CPU and GPU executions. The main computational kernels remain
resident in GPU memory, with element-level loops mapped to OpenACC gangs
and nodal operations mapped to vector-level parallelism. MPI-based domain
decomposition is used for multi-GPU execution.
The GPU implementation includes standard and split-form DGSEM volume
operators, multiple Riemann solvers, the BR1 viscous discretisation,
explicit time integrators, positivity-preserving limiting, LES models and
wall modelling. CPU and GPU implementations retain the same numerical
formulation and are checked through common verification and validation
cases.
The current GPU implementation and performance assessment primarily target
NVIDIA GPU systems through the NVIDIA HPC SDK OpenACC compiler. The solver
currently uses conforming, curvilinear hexahedral meshes, with mesh formats
including GMSH, HDF5 and SpecMesh/HOHQMesh.
\end{small}

\section{Introduction}

In the current era of artificial intelligence (AI), the quality and quantity of training data directly dictate the accuracy and generalisability of machine learning models. Although AI excels at pattern recognition and rapid inference, it remains fundamentally reliant on large volumes of high‐fidelity data that, for many fluid dynamic phenomena, can only be generated by accurate and efficient computational fluid dynamics (CFD) solvers~\cite{brunton2020machine}. 
Traditional low‐order, CPU‐bounded CFD codes struggle to keep up with these demands, leading to prohibitive compute times and limited parametric coverage.

CFD has become an indispensable tool across a wide range of engineering and scientific disciplines, underpinning the design and analysis of applications ranging from aerospace vehicles and wind turbines to cardiovascular flows. As the complexity and fidelity requirements of these simulations continue to grow, the computational cost of state‑of‑the‑art CFD solvers has risen dramatically. The growing role of CFD in early-stage design optimisation and uncertainty quantification, where ensembles of simulations must be performed within tight time frames, is becoming standard practice. In fields such as aerodynamic shape optimisation or wind farm layout design, the ability to complete hundreds or thousands of high‑resolution runs in parallel can mean the difference between timely, accurate answers and prohibitive compute costs. 

In parallel, high‑performance computing (HPC) architectures have evolved towards heterogeneous accelerator‐-based systems in which multiple Graphics Processing Units (multi‑GPU) nodes now dominate the Top500 and Green500 lists~\cite{top500_overview}. This convergence of escalating CFD demands and GPU‑centric hardware presents both a challenge and an opportunity: without careful porting and optimisation, traditional CPU‑centric solvers will fail to harness the full potential of modern supercomputers; conversely, GPU‑accelerated implementations can deliver orders‑of‑magnitude speedups, enabling simulations of unprecedented scale and fidelity.

Graphics Processing Units excel at data‑parallel workloads due to their massive on‑chip parallelism, high memory bandwidth, and specialised tensor and matrix engines. Yet fully exploiting these capabilities for CFD requires more than recompiling existing code; it demands rethinking data layouts, communication patterns, and algorithmic kernels to minimise memory traffic, overlap computation and data transfer, and adapt to the SIMD/SIMT (Single Instruction, Multiple Data/Multiple Threads) execution model. In addition, multi‑GPU clusters impose additional constraints on domain decomposition, halo‑exchange, and load balancing, especially for highly adaptive or unstructured discretisations. As CFD codes transition from single‑GPU demonstrations to production‑scale, multi‑GPU deployments, developers must address issues of scalability, portability, and maintainability to ensure robust performance on current and future exascale platforms.


High‐order discretisations for CFD --- such as discontinuous Galerkin (DG), spectral element methods (SEM), and flux‐reconstruction (FR) methods --- offer a compelling solution: they deliver superior per‐degree‐of‐freedom accuracy and reduced numerical dissipation, enabling coarser meshes (and thus fewer degrees of freedom) for a given solution fidelity~\cite{wang2013high}. 
However, their increased arithmetic intensity and more complex data‐access patterns also impose greater demands on memory bandwidth. 
Modern GPUs, with their thousands of computing cores and high‐bandwidth on‐chip memory, are ideally suited to element‐local, compute‐dominated kernels of high‐order schemes—but only if the solver is carefully ported and tuned for the GPU execution model.

Porting high-‐order CFD solvers to multi‐GPU systems is therefore not merely a performance exercise, but a strategic imperative in the AI‐driven age. Efficient GPU implementations can produce large volumes of high‐fidelity flow data at rates orders of magnitude faster than their CPU counterparts only, powering deep learning framework training and enabling real‐time digital twins. However, achieving this performance on scale requires overcoming challenges in data layout, kernel fusion, inter‐GPU communication, and load balancing, particularly for unstructured meshes and adaptive refinement. 

In this paper, we port the high order DG solver \texttt{HORSES3D} \cite{HORSES3D_paper} to multi-GPUs systems demonstrating how a high order finite element type solver can be restructured to new architectures using directive‐based offload (OpenACC) and domain decomposition techniques, and we report strong and weak scaling on up to 2048 NVIDIA H100 GPUs for problems exceeding ten billion degrees of freedom. Our results illustrate the path toward rapid large-ensemble CFD data generation to fuel the next wave of AI‐enabled discoveries in fluid dynamics.

\subsection{State of the art - high order solvers in GPUs}

High-order numerical methods—such as discontinuous Galerkin (DG), spectral element methods (SEM), and flux reconstruction (FR)—are particularly well suited for GPU-based acceleration due to their high arithmetic intensity and element-local computations. These features align well with GPU architectures, leading to substantial performance gains. A wide range of solvers have been developed to exploit this potential. For example, GALÆXI~\cite{Kurz2025}, an extension of FLEXI~\cite{krais2021flexi}, implements high-order DGSEM on GPUs and demonstrates excellent strong scaling up to 1024 GPUs, with more than $10^6$ degrees of freedom (DOFs) per device. NekRS~\cite{Fischer2023}, the GPUs-focused successor of Nek5000, adopts spectral elements and has been shown to scale to over 27,000 GPUs in Summit. PyFR~\cite{Vermeire2017} and ZEFR~\cite{Romero2020} implement FR schemes using OpenCL and CUDA, targeting compressible viscous flows.
Trixi.jl~\cite{schlottkelakemper2025trixi, ranocha2022adaptive, schlottkelakemper2021purely} implements high-order DG methods for conservation laws in the Julia programming language \cite{bezanson2017julia} and offers limited GPU support via KernelAbstractions.jl~\cite{Churavy_KernelAbstractions_jl}.
In contrast, other DGSEM codes like Nektar++~\cite{moxey2020nektar,cantwell2015nektar} and FLEXI~\cite{krais2021flexi} have remained CPU-oriented.
These codes use a variety of programming frameworks and strategies to exploit GPUs efficiently. CUDA and HIP remain the most commonly used for high performance, but directive-based approaches (e.g. OpenACC) and portability layers (e.g. Kokkos~\cite{trott2021kokkos}, RAJA~\cite{beckingsale2019raja}, SYCL~\cite{alpay2020sycl}) are also widely used. NekRS relies on OCCA to generate portable kernels for both CUDA and HIP targets, while Gasparino \emph{et al.}~\cite{Gasparino2024} used OpenACC to accelerate SOD2D, an SEM solver with minimal changes to the original code. In all cases, high efficiency is achieved by optimising memory access patterns (e.g., one element per thread block), fusing kernels, and overlapping computation and communication.

Performance results consistently highlight the advantages of GPU acceleration. Speedups of over $3\times$ compared to dual-socket CPU nodes have been reported on single GPUs like V100~\cite{Gasparino2024}, and newer architectures (e.g., A100, MI250X) provide further gains—up to $1.6\times$ in throughput for the same problem size. 

The range of applications of GPU-accelerated high-order solvers is broad and continues to expand. In aerodynamics and turbomachinery, these methods are increasingly used for complex compressible-flow configurations, as illustrated by GALÆXI simulations of the NASA Rotor 37 compressor rotor~\cite{Kurz2025}. PyFR has also demonstrated the viability of GPU-accelerated high-order flux-reconstruction methods for scale-resolving simulations of turbulent flows over wings, airfoils, and control-surface configurations~\cite{witherden2025pyfr,Romero2020}. Beyond external aerodynamics, canonical turbulence benchmarks, including the Taylor--Green vortex, turbulent jets, and channel-flow configurations, have been used extensively to assess the accuracy, robustness, and scalability of high-order GPU-enabled solvers such as SOD2D and related DG/FR frameworks~\cite{Gasparino2024,Bull2015,Karakus2019}. 
Across all these domains, GPU acceleration is enabling simulations with finer resolution and longer time horizons than were previously feasible, dramatically expanding the capabilities of high-order computational fluid dynamics.

\section{\texttt{HORSES3D} on multi-GPU systems}

\subsection{\texttt{HORSES3D}: High order solver}
\texttt{HORSES3D} is an open source solver 
\cite{HORSES3D_paper}
developed at the ETSIAE-UPM School of Aeronautics in Madrid, and available on Github (\url{https://github.com/horses-framework/HORSES3D}) with the GPU ported version available in (\url{https://github.com/horses-framework/HORSES3D-gpu}).
\texttt{HORSES3D} is a high-order discontinuous Galerkin (DG) framework with capabilities covering a wide range of flow regimes and physical models, including incompressible flows~\cite{manzanero2020entropy}, compressible and supersonic flows~\cite{ntoukas2025comparative,mateo2025unsupervised, lodares2022entropy, rueda2019p}, multiphase flows~\cite{ballout2025acoustic, ntoukas2022entropy}, and aeroacoustics~\cite{BOTEROBOLIVAR2024120476,oscarAL}. In the present work, we focus on the compressible-flow formulation to assess and illustrate the performance of the GPU implementation. Note that the incompressible and multiphase formulations have also been ported to GPU architectures, although their detailed performance analysis is beyond the scope of this paper.


\texttt{HORSES3D} is an h/p discontinuous Galerkin solver designed for high-order simulations of fluid dynamics. The "h" in h/p refers to mesh refinement (h-refinement), while "p" refers to increasing the polynomial order of the solution approximation (p-refinement). This dual approach allows \texttt{HORSES3D} to achieve high accuracy with fewer degrees of freedom compared to traditional low-order methods, making it particularly effective for capturing the fine details of the complex, multi-scale flow phenomena (e.g., turbulent flows) \cite{kessasra2024comparison}. One of the main advantages of DG methods is their ability to accurately capture high-order spatial and temporal variations of the solution, which makes them particularly suitable for simulating flows with sharp gradients and complex flow phenomena. DG methods also exhibit good numerical stability and conservation properties because of the local nature of the approximation and the use of fluxes at the interfaces of the elements. 

%

Here, we provide only a brief overview of the fundamental concepts of DG discretisations for the compressible Navier--Stokes (NS) retained in this work; see the Appendix \ref{sec:cNS}.
The physical domain is subdivided into non-overlapping curvilinear hexahedral elements, $e$, which are geometrically transformed to a reference element, $el$, using a polynomial transfinite mapping that relates the physical coordinates $\vec{x}$ and the local reference coordinates $\vec{\xi}$. This transformation is applied to Eq.~\eqref{eq:compressibleNScompact_transformed}, resulting in the following:

\begin{equation}
\boldsymbol{q}_t  + \nabla_x\cdot\ssvec{F}_e = \nabla_x\cdot\ssvec{F}_{v,turb}+\boldsymbol{S({q})},
\rightarrow
J \boldsymbol{q}_t  + \nabla_\xi\cdot\cssvec{F}_e = \nabla_\xi\cdot\cssvec{F}_{v,turb}+J \boldsymbol{S({q})},
\label{eq:compressibleNScompact_transformed}
\end{equation}
where $J$ is the Jacobian determinant of the transfinite mapping, $\nabla_\xi$ is the differential operator in the reference space and $\cssvec{F}$ are the contravariant fluxes \cite{kopriva2009implementing}. 

To derive DG schemes, we multiply Eq.~\eqref{eq:compressibleNScompact_transformed} by a locally smooth test function $\phi$. We integrate over an element $el$ to obtain the weak form:
\begin{equation}\label{eq::NS2}
\int_{el}J \boldsymbol{q}_t\phi \, \d \vec{x}
+\int_{el} \nabla_\xi\cdot\cssvec{F}\phi \, \d \vec{x}  =
\int_{el}J\boldsymbol{S({q})}\phi \, \d \vec{x},
\end{equation}
where we have gathered the contributions of inviscid and viscous terms in the flux $\cssvec{F} := \cssvec{F}_e - \cssvec{F}_{v,turb}$. By integrating the flux term by parts, we obtain a local weak form of the equations (one per mesh element) with the boundary fluxes separated from the interior,
\begin{equation}\label{eq::NS3}
\int_{el}J \boldsymbol{q}_t\phi \, \d \vec{x} +  
\int_{\partial el} \cssvec{F}\cdot\hat{\mathbf{n}}\phi \, \d \vec{s}
-\int_{el} \cssvec{F}\cdot\nabla_\xi\phi \, \d \vec{x}
=
\int_{el}J\boldsymbol{S({q})}\phi \, \d \vec{x},
\end{equation}
where $\hat{\mathbf{n}}$ is the unit outward-pointing vector of each face of the reference element ${\partial el}$.
Since the solution is allowed to be discontinuous at inter-element faces, surface fluxes are replaced with a numerical flux $\mathbf{F}_{e}^{\star}$, to couple the elements:
\begin{equation}\label{eq::NS4}
\int_{el}J \boldsymbol{q}_t\phi \, \d \vec{x} +  
\int_{\partial el} \cssvec{F}^{\star}\cdot\hat{\mathbf{n}}\phi \, \d \vec{s}
-\int_{el} \cssvec{F}\cdot\nabla_\xi\phi \, \d \vec{x}
=
\int_{el}J\boldsymbol{S({q})}\phi \, \d \vec{x},
\end{equation}

The equations for each element are coupled with the equations of the neighboring elements through the numerical fluxes, which are composed of inviscid and viscous parts, $\mathbf{F}^{\star} := \mathbf{F}^{\star}_e - \mathbf{F}^{\star}_{v,turb}$.
The inviscid numerical flux function, $\mathbf{F}^{\star}_e$, is computed with an approximate Riemann solver (see, e.g., \cite{toro2013riemann}), and the viscous numerical flux function, $\mathbf{F}^{\star}_{v,turb}$, is computed using a consistent discontinuous Galerkin discretization \cite{arnold2001dg,HORSES3D_paper}. 
Non-linear inviscid and viscous numerical fluxes (including turbulent ones) can be chosen appropriately to control dissipation in the numerical scheme~\cite{Manzanero_2020,Ferrer_2017,jumpKou}.

A common strategy consists in using integration by parts on the flux volume integral again to arrive at the so-called weak-strong form of the DG discretisation:
\begin{equation}\label{eq::NS5}
\int_{el}J \boldsymbol{q}_t\phi \, \d \vec{x} +  
\int_{\partial el} \left(\cssvec{F}^{\star} - \cssvec{F}\right) \cdot\hat{\mathbf{n}}\phi \, \d \vec{s}
+\int_{el} \nabla_{\xi}\cdot\cssvec{F}\phi \, \d \vec{x}
=
\int_{el}J\boldsymbol{S({q})}\phi \, \d \vec{x}.
\end{equation}

In a final step, we approximate the numerical solution and fluxes using polynomials (of order $p$) and evaluate all integrals using Gaussian quadrature rules.
In \texttt{HORSES3D-GPU}, we use a collocated discontinuous Galerkin spectral element method (DGSEM) \cite{black1999conservative, kopriva2009implementing} and can select either Gauss--Legendre or Gauss--Lobatto--Legendre collocation points.
For the latter choice of collocation points, we have implemented the so-called split-form (or flux-differencing) form of the differentiation operator \cite{gassner2016split, carpenter2013high, carpenter2013high} based on two-point fluxes, which allows the selection of several kinetic-energy-preserving and entropy-conservative volume numerical fluxes. The implementation of the split-form DGSEM in \texttt{HORSES3D-GPU} includes several performance optimizations, which are detailed in \cite{ranocha2023efficient}.

Several choices of split forms and Riemann fluxes are included and have been ported to the GPU version of the code; see Section \ref{sec:gpu_port}.  In addition, we have complemented the compressible Navier--Stokes equations with the Vreman \cite{Vreman_2004} and Wale \cite{nicoud1999subgrid} Large Eddy Simulation subgrid models. Note that the source term $\mathbf{S}$ can be used to incorporate additional physics, such as immersed boundaries (mesh free methods) \cite{kou2022immersed,HORSES3D_paper,colombo2025high} or actuator lines to simulate wind turbines ~\cite{Marino2024,Botero-Bolivar2024}. Finally, to advance the solution in time, we use explicit low-storage Runge-Kutta time marching schemes.

\subsection{Preliminaries: Selecting a GPU paradigm}
There are several approaches for porting scientific codes to GPU architectures, ranging from low-level, vendor-specific programming models to higher-level directive-based or performance-portability frameworks. These approaches differ in terms of development effort, required code restructuring, degree of control over the GPU memory hierarchy and execution model, achievable performance, and long-term maintainability. At the lowest level, CUDA remains the native programming model for NVIDIA GPUs~\cite{nvidia_cuda_guide}, while HIP provides a CUDA-like C++ programming interface within the AMD ROCm ecosystem and is commonly used to target AMD GPUs while facilitating the porting of CUDA codes~\cite{amd_rocm_programming_guide,amd_hip_programming_model}. Higher-level alternatives include directive-based models such as OpenACC and OpenMP target offloading~\cite{openacc_guide,openmp_api}, as well as performance-portability frameworks such as Kokkos, RAJA, and SYCL, which aim to reduce vendor lock-in by enabling a single source code to target multiple CPU and GPU back ends~\cite{edwards2014kokkos,davis2024gpuportability,khronos_sycl}.

Low-level GPU programming models, such as CUDA and HIP, extend popular languages such as C/C++ and Fortran and provide fine-grained control over the GPU execution model, memory hierarchy, and hardware resources through dedicated APIs. This approach can deliver the highest performance when implemented by expert GPU developers, but it often requires substantial code restructuring and may lead to the parallel maintenance of multiple code paths for CPUs and different GPU vendors. It is therefore particularly attractive for relatively mature codebases whose numerical kernels are stable and where peak performance is the dominant requirement. 
At a higher level, performance-portability approaches such as OpenCL, SYCL, Kokkos, and RAJA aim to reduce vendor lock-in by enabling a single source code to target different CPU and GPU back ends. These frameworks generally require less hardware-specific programming than CUDA or HIP, while still allowing good control over parallel execution and data management. As a result, they offer an attractive compromise between performance, portability, and maintainability, especially for large scientific codes intended to run on heterogeneous HPC platforms. However, most of these approaches are primarily designed around C and C++ programming models, and their support for Fortran codes is either limited, indirect, or unavailable, which can be a significant constraint for legacy scientific software.

\subsection{Porting \texttt{HORSES3D} to GPUs using openACC}

The porting process of \texttt{HORSES3D} to GPUs had some requirements that led the development team to select OpenACC. 
Since multiple developers work on this multiphysics CFD framework, the GPU version should enable a quick and straightforward process. 
OpenACC works through directives on top of the original code, so development time is kept to a minimum, and there is no need for every developer to be an expert in CUDA/HIP to create an equivalent GPU version.
The back-end should support the already existing data structure of the code with little modifications. OpenACC is able to handle complex data structures and data stored in AoS (array of structures) format. This was the deciding factor for selecting OpenACC over OpenMP, as OpenMP offloading to GPUs does not offer support for data structures as complex as in \texttt{HORSES3D}.

In the \texttt{HORSES3D} GPU solver, all data are initialised on the CPU and then transferred to the GPU before starting the time-stepping procedure. As in other high-performance GPU-resident CFD frameworks, all computationally intensive operations have been ported to the GPU to avoid the bottleneck associated with slow CPU--GPU data transfers. An example of the ported code with OpenACC is presented in Listings \ref{lst:code_outer} and \ref{lst:code_inner_std}. The element and face level loops are exposed to the gang level of parallelisation, where each gang is associated with one element. A gang corresponds to a thread-block in CUDA/HIP terminology. This is presented in Listing \ref{lst:code_outer} for a volumetric integral. Then for the inner loop, presented in Listing \ref{lst:code_inner_std}, we use vector loops to associate the nodal points within each element with each thread (compute unit). Through this simple approach, we have been able to effectively port the \texttt{HORSES3D} code to GPUs and exploit the significant performance advantage that modern GPU hardware offers. 
The use of the variable \texttt{r\_volInt} reduces the successive reads and writes to global memory and increases the code efficiency.
We have used this strategy in several parts of \texttt{HORSES3D} and re-designed the loops to improve the memory management. 
This has been necessary to extract the maximum available performance from \texttt{HORSES3D} on GPUs and have a framework that is easy to maintain and has performance comparable to other known CFD frameworks that have native CUDA/HIP support. 

More generally, the porting effort prioritises performance-critical kernels, including the evaluation of volume integrals, surface fluxes, viscous terms, and time integration. These kernels are executed on the GPU with data that reside entirely in device memory, minimising data transfers between the host and the device.

\begin{listing}[H]
\begin{lstlisting}[caption={Example of outer loop parallelisation with OpenACC.}, label={lst:code_outer}]
subroutine TimeDerivative_VolumetricContribution(mesh)

    !...Operations before...

    !$acc parallel loop gang vector_length(128) present(mesh) async(1)
    do eID = 1 , size(mesh % elements)

        !...Operations before...
            
        call ScalarWeakIntegrals_StdVolumeGreen(mesh % elements(eID) % Nxyz, NCONS, &
                                                mesh % elements(eID) % storage % contravariantFlux, &
                                                mesh % elements(eID) % storage % QDot)
    end do
    !$acc end parallel loop 

end subroutine TimeDerivative_VolumetricContribution
\end{lstlisting}
\end{listing}

\begin{lstlisting}[caption={Example of inner loop parallelisation with OpenACC: Standard weak divergence operator.}, label={lst:code_inner_std}]
subroutine ScalarWeakIntegrals_StdVolumeGreen( Nxyz, NEQ, F, volInt )
   !$acc routine vector
   implicit none
   integer,             intent(in)  :: Nxyz(3)
   integer,             intent(in)  :: NEQ
   real(kind=RP),       intent(in)  :: F     (1:NCONS, 0:Nxyz(1), 0:Nxyz(2), 0:Nxyz(3), 1:NDIM)
   real(kind=RP),    intent(inout)  :: volInt(1:NCONS, 0:Nxyz(1), 0:Nxyz(2), 0:Nxyz(3))
!
!  ---------------
!  Local variables
!  ---------------
!
   integer   :: i, j, k, l, eq
   real(kind=RP) :: r_volInt

   !$acc loop vector collapse(4) private(r_volInt)
      do k = 0, Nxyz(3) ; do j = 0, Nxyz(2) ; do i = 0, Nxyz(1) ; do eq = 1, NEQ
      
      r_volInt = volInt(eq,i,j,k)
      
      !$acc loop seq 
      do l = 0, Nxyz(1)
         r_volInt = r_volInt + NodalStorage(Nxyz(1)) % hatD(i,l) * F(eq,l,j,k,IX)
      end do  
      
      !$acc loop seq 
      do l = 0, Nxyz(2)
         r_volInt = r_volInt + NodalStorage(Nxyz(2)) % hatD(j,l) * F(eq,i,l,k,IY)
      end do  
      
      !$acc loop seq 
      do l = 0, Nxyz(3)
         r_volInt = r_volInt + NodalStorage(Nxyz(3)) % hatD(k,l) * F(eq,i,j,l,IZ)
      end do             
      
      ! Write back to global memory
      volInt(eq,i,j,k) = r_volInt        
   end do ; end do ; end do ; end do

end subroutine ScalarWeakIntegrals_StdVolumeGreen
\end{lstlisting}

The split-form volume integral, presented in Listing~\ref{lst:code_inner_split}, follows the same gang/vector mapping as the standard weak-form kernel: one gang per element and one vector lane per nodal point $(i,j,k)$, with the summation index $l$ and the equation index \texttt{eq} kept as sequential loops within each lane. 
Unlike the standard divergence operator, the two-point flux entering this integral is evaluated on the fly for every pair of nodes along each coordinate line rather than read from a precomputed tensor, since materialising the full two-point flux array in global memory for every element would introduce a memory-bandwidth bottleneck that offsets the benefit of the split-form scheme. 
As in Listing \ref{lst:code_inner_std}, the partial sum over $l$ is accumulated in a register-resident variable before being written back to \texttt{QDot} once per node, avoiding repeated global-memory reads and writes at every step of the summation.

Two-point fluxes of this type are symmetric, $F(Q_i,Q_l) = F(Q_l,Q_i)$, a property that can be exploited on CPU to halve the number of flux evaluations by scattering each computed pair to both output indices \cite{ranocha2023efficient}.
This is deliberately not exploited on GPU: scattering a single evaluation to two output nodes would require two different vector lanes to update each other's contribution to \texttt{QDot} within the same vectorised loop, which is only safe with atomic memory operations. 
Since atomic throughput on current GPU hardware is substantially lower than that of private, per-lane accumulation, the redundant flux evaluations avoided by exploiting symmetry are, in this setting, cheaper than the synchronisation cost of exploiting it. 
\texttt{HORSES3D} therefore evaluates the two-point flux independently for both orderings of each node pair and keeps the accumulation strictly private to a single vector lane.
The advective split-form contribution and the viscous weak-form contribution are further fused into the accumulation performed within a single kernel launch, rather than being computed in two separate \texttt{!\$acc parallel} regions, which removes an additional full read-modify-write pass over \texttt{QDot} per element per Runge--Kutta stage.

\begin{lstlisting}[caption={Example of inner loop parallelisation with OpenACC: Split-form divergence operator.}, label={lst:code_inner_split}]
subroutine ScalarWeakIntegrals_SplitVolumeDivergence( e, Fv, QDot )
   !$acc routine vector
   use ElementClass
   use RiemannSolvers_NS
   implicit none
   type(Element), intent(in)    :: e
   real(kind=RP), intent(in)    :: Fv  (1:NCONS, 0:e%Nxyz(1), 0:e%Nxyz(2), 0:e%Nxyz(3), 1:NDIM)
   real(kind=RP), intent(inout) :: QDot(1:NCONS, 0:e%Nxyz(1), 0:e%Nxyz(2), 0:e%Nxyz(3))
!
!  ---------------
!  Local variables
!  ---------------
!
   integer       :: i, j, k, l, eq, Nx, Ny, Nz
   real(kind=RP) :: Fs(1:NCONS)
   real(kind=RP) :: r_QDot(1:NCONS)

   Nx = e % Nxyz(1) ; Ny = e % Nxyz(2) ; Nz = e % Nxyz(3)

   !$acc loop vector collapse(3) private(Fs, r_QDot)
   do k = 0, Nz ; do j = 0, Ny ; do i = 0, Nx

      !$acc loop seq
      do eq = 1, NCONS
         r_QDot(eq) = QDot(eq,i,j,k)          ! single read from global memory
      end do

      ! xi
      !$acc loop seq
      do l = 0, Nx
         call TwoPointFlux_Selector(e % storage % Q(:,i,j,k), e % storage % Q(:,l,j,k), &
                                    e % geom % jGradXi(:,i,j,k), e % geom % jGradXi(:,l,j,k), Fs )
         !$acc loop seq
         do eq = 1, NCONS
            r_QDot(eq) = r_QDot(eq) - NodalStorage(Nx) % sharpD(i,l) * Fs(eq) &
                                    + NodalStorage(Nx) % hatD(i,l)   * Fv(eq,l,j,k,IX)
         end do
      end do

      ! eta
      !$acc loop seq
      do l = 0, Ny
         call TwoPointFlux_Selector( e % storage % Q(:,i,j,k), e % storage % Q(:,i,l,k), &
                                    e % geom % jGradEta(:,i,j,k), e % geom % jGradEta(:,i,l,k), Fs )
         !$acc loop seq
         do eq = 1, NCONS
            r_QDot(eq) = r_QDot(eq) - NodalStorage(Ny) % sharpD(j,l) * Fs(eq) &
                                    + NodalStorage(Ny) % hatD(j,l)   * Fv(eq,i,l,k,IY)
         end do
      end do

      ! zeta
      !$acc loop seq
      do l = 0, Nz
         call TwoPointFlux_Selector( e % storage % Q(:,i,j,k), e % storage % Q(:,i,j,l), &
                                    e % geom % jGradZeta(:,i,j,k), e % geom % jGradZeta(:,i,j,l), Fs )
         !$acc loop seq
         do eq = 1, NCONS
            r_QDot(eq) = r_QDot(eq) - NodalStorage(Nz) % sharpD(k,l) * Fs(eq) &
                                    + NodalStorage(Nz) % hatD(k,l)   * Fv(eq,i,j,l,IZ)
         end do
      end do

      !$acc loop seq
      do eq = 1, NCONS
         QDot(eq,i,j,k) = r_QDot(eq)          ! single write back to global memory
      end do

   end do ; end do ; end do

end subroutine ScalarWeakIntegrals_SplitVolumeDivergence
\end{lstlisting}

\subsection{Overall GPU Implementation and Ported Functionalities}
\label{sec:gpu_port}

The GPU version of \texttt{HORSES3D} preserves the numerical formulation of the original CPU solver while restructuring the computational kernels for efficient execution on GPU architectures. The porting strategy relies on directive-based parallelisation to maintain a single code base, ensuring that the physical models and numerical discretisations remain identical in CPU and GPU executions. All features listed in Table~\ref{tab:gpu_port} have been ported to GPUs and verified through dedicated test cases included in the \texttt{HORSES3D} repository. These tests are part of the continuous-integration workflow and can be executed on both CPU and GPU backends, enabling direct numerical-consistency checks across architectures.

The configuration considered in the present work represents a subset of these capabilities. In particular, the Taylor--Green vortex and turbulent channel-flow cases and the full aircraft (Common Research Model) rely on combinations of split forms, Riemann solvers, and turbulence models that are representative of production CFD simulations, while remaining fully supported and validated in the GPU implementation.

\begin{table}[h!]
\centering
\caption{Functionalities ported to the GPU version of \texttt{HORSES3D} and validated through CPU/GPU cases.}
\label{tab:gpu_port}

\begin{tabular}{ll}
\toprule
\textbf{Category} & \textbf{Implemented features (GPU)} \\
\midrule

\textbf{Choice of nodes} &
Gauss--Legendre (GL), Gauss--Legendre--Lobatto (GLL) \\

\midrule

\textbf{Split forms}
&
Standard~\cite{10.1016/j.jcp.2016.09.013},
Morinishi~\cite{morinishi2010skew,10.1016/j.jcp.2016.09.013}, \\
(for GLL)
&
Ducros~\cite{DUCROS2000114,10.1016/j.jcp.2016.09.013},
Kennedy--Gruber~\cite{kennedy2008reduced}, \\
&
Pirozzoli~\cite{pirozzoli2010generalized},
Ismail--Roe~\cite{ismail2009affordable}, 
Chandrasekhar~\cite{chandrashekar2013kinetic} \\

\midrule

\textbf{Riemann solvers}
&
Central (non-dissipative), Roe~\cite{toro2013riemann}, \\
&
Low-dissipation Roe~\cite{osswald2016l2roe},
Roe--Pike~\cite{roe1984}, \\
&
Matrix dissipation~\cite{ismail2009affordable},
Local Lax--Friedrichs/Rusanov~\cite{Friedrichs1971}\\
\midrule

\textbf{Viscous discretisation} &
Bassi--Rebay 1 (BR1)~\cite{bassi1997high} \\

\midrule

\textbf{Time integration}
&
Explicit Euler, 
Low-storage RK3~\cite{williamson1980low}, \\
&
Low-storage RK5 (Carpenter--Kennedy)~\cite{carpenter1994rk}, \\
&
SSPRK33~\cite{shu1988}, 
SSPRK43~\cite{spiteri2002} \\

\midrule

\textbf{Limiters} &
Positivity-preserving limiter~\cite{zhang2011} \\

\midrule

\textbf{Turbulence models} &
WALE~\cite{nicoud1999subgrid},
Vreman~\cite{vreman2004eddy} \\

\midrule

\textbf{Wall modelling} &
Reichardt wall model~\cite{frere2017application} \\

\bottomrule
\end{tabular}
\end{table}


\section{Verification and validation}

The GPU implementation of \texttt{HORSES3D} retains the numerical implementation of the original CPU code\cite{HORSES3D_paper}. The porting process consists of restructuring and annotating the existing computational kernels for execution on GPU devices, while maintaining the definitions for the spatial discretisation, time-integration schemes, numerical fluxes, boundary conditions, and physical models. Therefore, the objective of this section is not to revalidate the full numerical methodology of \texttt{HORSES3D}, which has already been extensively assessed in previous work \cite{ferrer2023high}, but rather to verify that the GPU-accelerated implementation preserves the accuracy and physical behaviour of the original solver.

For completeness, three representative verification and validation cases are included for the 3D compressible Navier--Stokes equations. 
First, a method of manufactured solutions (MMS) test is used to verify the expected high-order convergence of the GPU implementation. Second, the Taylor--Green vortex at $\mathrm{Re}=1600$ is considered to validate the solver against established DNS/reference data for transitional and turbulent flows. Additionally, a turbulent channel-flow including wall modelling is used to assess the behaviour of the GPU solver in a wall-bounded turbulent-flow case. 


\subsection{Verification: Method of Manufactured Solutions}
\label{sec:MMS}

To verify that the GPU-accelerated implementation of \texttt{HORSES3D} preserves the expected order of accuracy of the underlying DG discretisation, a method of manufactured solutions (MMS) test is performed. Since the GPU port retains the same spatial discretisation, time-integration schemes, and numerical fluxes as the CPU version, this test focuses on confirming that the GPU-specific optimisations do not alter the numerical convergence properties of the solver.

A smooth analytical solution is prescribed for all flow variables, and the corresponding source term is added to the compressible Navier--Stokes equations such that the manufactured solution satisfies the modified governing equations exactly. The analytical expression for the manufactured solution is defined as:

\begin{equation}
\begin{aligned}
    \rho (x,y,z,t) &= 1.0, \\
    \rho u (x,y,z,t) &= \sin(t) \sin(\pi x) \cos(\pi y) \cos(\pi z), \\
    \rho v (x,y,z,t) &= \sin(t) \sin(\pi y) \cos(\pi x) \cos(\pi z), \\
    \rho w (x,y,z,t) &= -2 \sin(t) \sin(\pi z) \cos(\pi x) \cos(\pi y), \\
    \rho e (x,y,z,t) &= 0.5 \sin^2(t) \sin^2(\pi x) \cos^2(\pi y) \cos^2(\pi z) + 0.5 \sin^2(t) \sin^2(\pi y) \cos^2(\pi x) \cos^2(\pi z) + \\ & \quad \, \, 2 \sin^2(t) \sin^2(\pi z) \cos^2(\pi x) \cos^2(\pi y) + 0.25 \cos(t) \sin(\pi x) \sin(\pi y) \sin(\pi z) + 2.5,
    \end{aligned}
    \label{eq:MMS_analytical_solution}
\end{equation}

with $\rho$ as the density, $(u,v,w)$ as the three components of the velocity, and $e$ as the total energy.
The associated forcing terms are omitted here for brevity.

\paragraph{Numerical setup:}
The MMS test is carried out on a sequence of uniform Cartesian meshes with hexahedral elements $3^3$, $4^3$, $5^3$, and $6^3$, within a domain defined as a cube of size $L=2$ with periodic boundary conditions on all faces. For each mesh, simulations are performed using increasing polynomial orders, and the numerical solution is advanced until convergence. The $L^2$ error norm of the conservative variables is then computed with respect to the manufactured analytical solution.

The compressible Navier--Stokes equations are solved, for a Mach number of $M=0.3$ and a Reynolds number of $Re=100000$, using the DGSEM formulation with Gauss--Lobatto nodes and the Roe Riemann solver for inter-element fluxes. Viscous terms are discretised using the BR1 scheme. Time integration is performed using an explicit third-order Runge--Kutta method with a fixed time step $\Delta t = 10^{-8}$.

\begin{figure}
    \centering
    \includegraphics[width=0.9\linewidth]{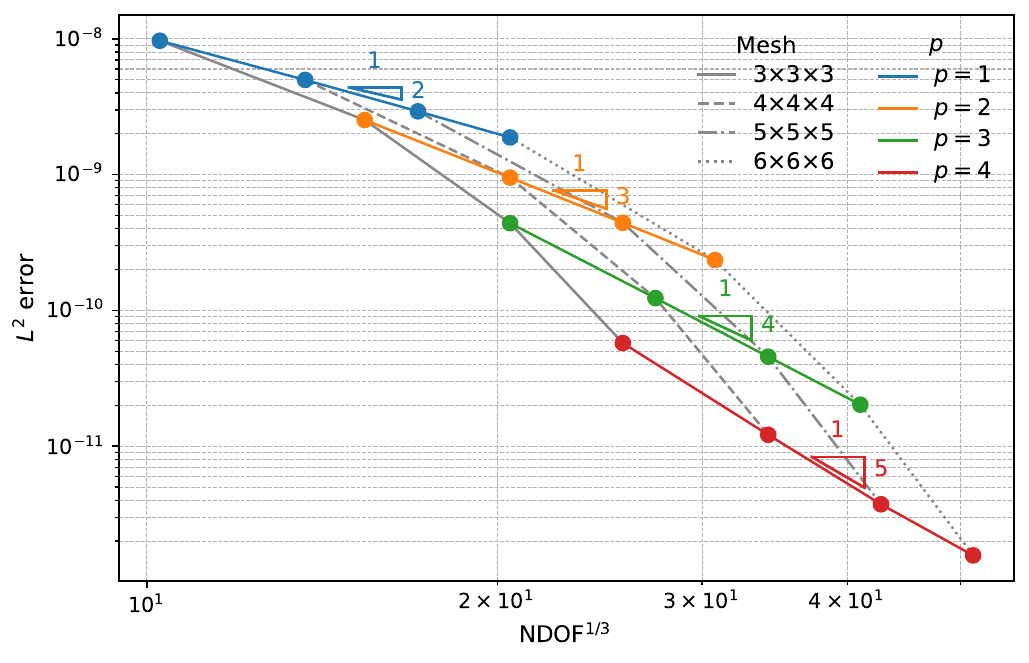}
    \caption{Convergence of the $L^2$ error against the number of degrees of freedom for different meshes and polynomial orders in the manufactured solution test case.} 
    \label{fig:convergence_MS}
\end{figure}

\paragraph{Results:} Figure~\ref{fig:convergence_MS} shows the exponential convergence, for h- and p-refinement, characteristic of high order solvers. The figure shows the error $L^2$ as a function of the square root of the total number of degrees of freedom for the different meshes and polynomial orders. For each fixed mesh, increasing the degree of the polynomial leads to an exponential reduction in the error, whereas mesh refinement for a fixed polynomial order yields consistent error decay. The reported convergence rates are in agreement with the theoretical expectations for high order formulations, confirming that the GPU-accelerated solver achieves the designed high order accuracy. Furthermore, the validation test case has been run both in CPU and GPU, providing the same results up to double-precision accuracy.
These results demonstrate that the GPU implementation of \texttt{HORSES3D} reproduces the same convergence behaviour as the original CPU solver. 

\subsection{Validation: Taylor--Green Vortex at $\mathrm{Re}=1600$}
\label{sec:validation_TGV}

The three-dimensional Taylor--Green vortex (TGV) at Reynolds number $\mathrm{Re}=1600$ is used to validate the GPU-accelerated version of \texttt{HORSES3D} to compute transitional and turbulent flows. This test case is widely adopted as a benchmark for assessing the ability of numerical schemes to reproduce the laminar--turbulent transition and the subsequent development of homogeneous isotropic turbulence.

The initial condition is given by the standard divergence-free TGV velocity field,
\begin{equation}
\begin{aligned}
\rho &= \rho_0, \\
u &= u_0 \sin(x/L_0)\cos(y/L_0)\cos(z/L_0), \\
v &= -u_0 \cos(x/L_0)\sin(y/L_0)\cos(z/L_0), \\
w &= 0,
\end{aligned}
\end{equation}
with the corresponding pressure field
\begin{equation}
p =
\frac{\rho_0 u_0^2}{\gamma M_0^2}
+
\frac{\rho_0 u_0^2}{16}
\left[
\left(
\cos(2x/L_0)+\cos(2y/L_0)
\right)
\left(
\cos(2z/L_0)+2
\right)
\right],
\end{equation}
where $\gamma=1.4$, $M_0=0.1$, $\rho_0=1$, $u_0=1$, and $L_0=1$. Time is non-dimensionalised using the convective time scale $t_c=L_0/u_0$.\\

\paragraph{Numerical setup:}
The computational domain is a periodic cube of size $L=2\pi$, discretised using $16^3$ uniform hexahedral elements with polynomial order $P=7$, corresponding to approximately $2.1\times10^6$  degrees of freedom. The compressible Navier--Stokes equations are solved using the DGSEM formulation with Gauss--Lobatto nodes, the Chandrasekhar split form, and the Lax-Friedrichs  solver for inter-element fluxes. Viscous terms are discretised using the BR1 scheme. Time integration is performed using an explicit third-order Runge--Kutta method with a fixed non-dimensional time step $\Delta t = 2\times10^{-4}$. No explicit subgrid-scale model is used, resulting in an implicit LES configuration. This configuration is considered here because it provides a representative baseline set of coefficients for \texttt{HORSES3D}. The full set of combinations of split forms and Riemann fluxes available in the GPU-ported version of \texttt{HORSES3D} has been systematically analysed in a separate study~\cite{rubio2025can}.

We verify the accuracy of the solver by computing the kinetic energy as
\begin{equation}
E =
\frac{1}{\rho_0 V}
\int_V
\frac{1}{2}\rho \mathbf{u}\cdot\mathbf{u}\,dV,
\end{equation}
and the kinetic energy dissipation rate, which is defined as $-\mathrm{d}E/\mathrm{d}t$.

\begin{figure}[h!]
    \centering
    \begin{subfigure}[t]{0.49\textwidth}
        \centering
        \includegraphics[width=\textwidth]{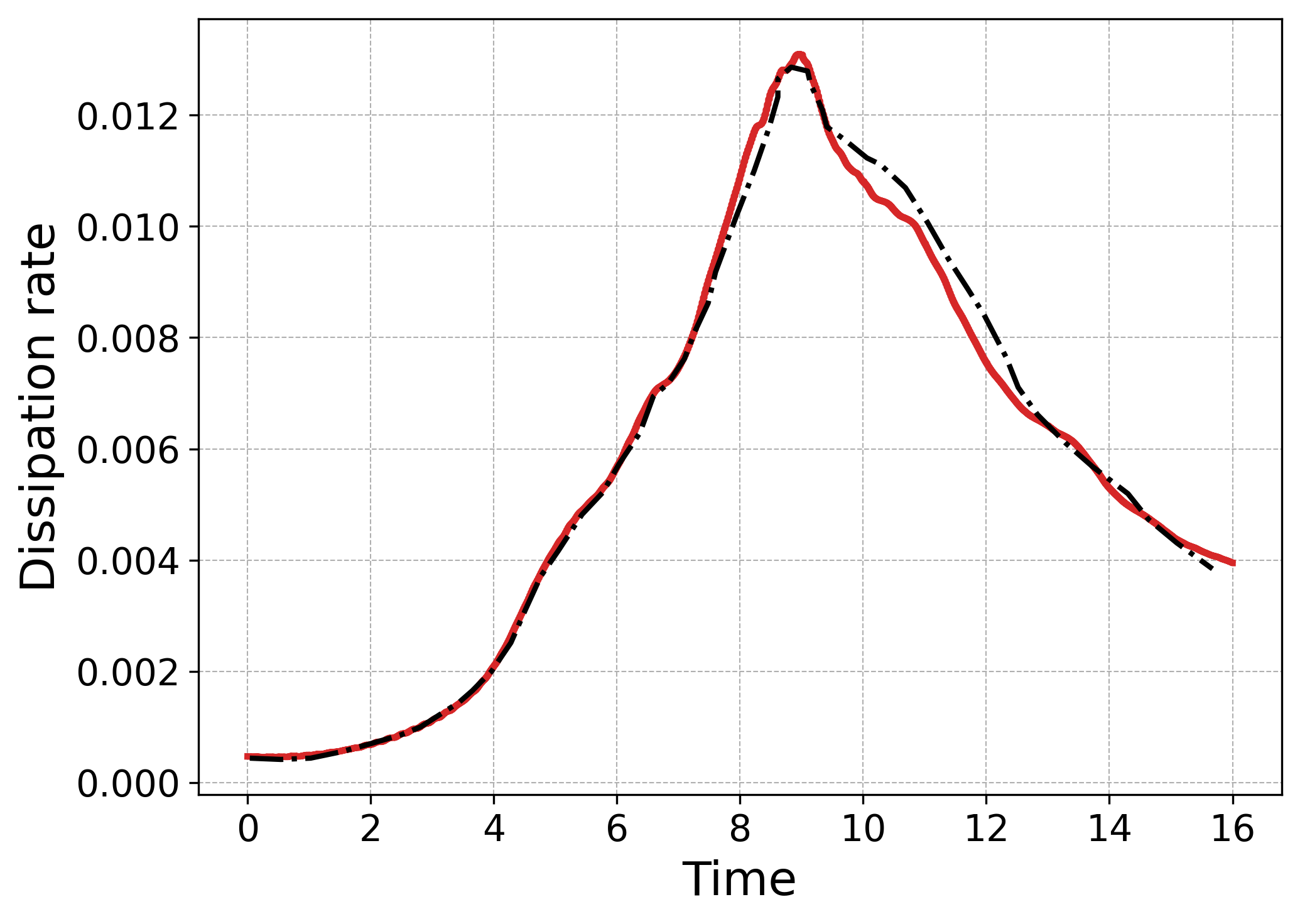}
        \caption{Kinetic energy dissipation rate.}
        \label{fig:TGV1600_dissipation}
    \end{subfigure}
    \hfill
    \begin{subfigure}[t]{0.49\textwidth}
        \centering
        \includegraphics[width=\textwidth]{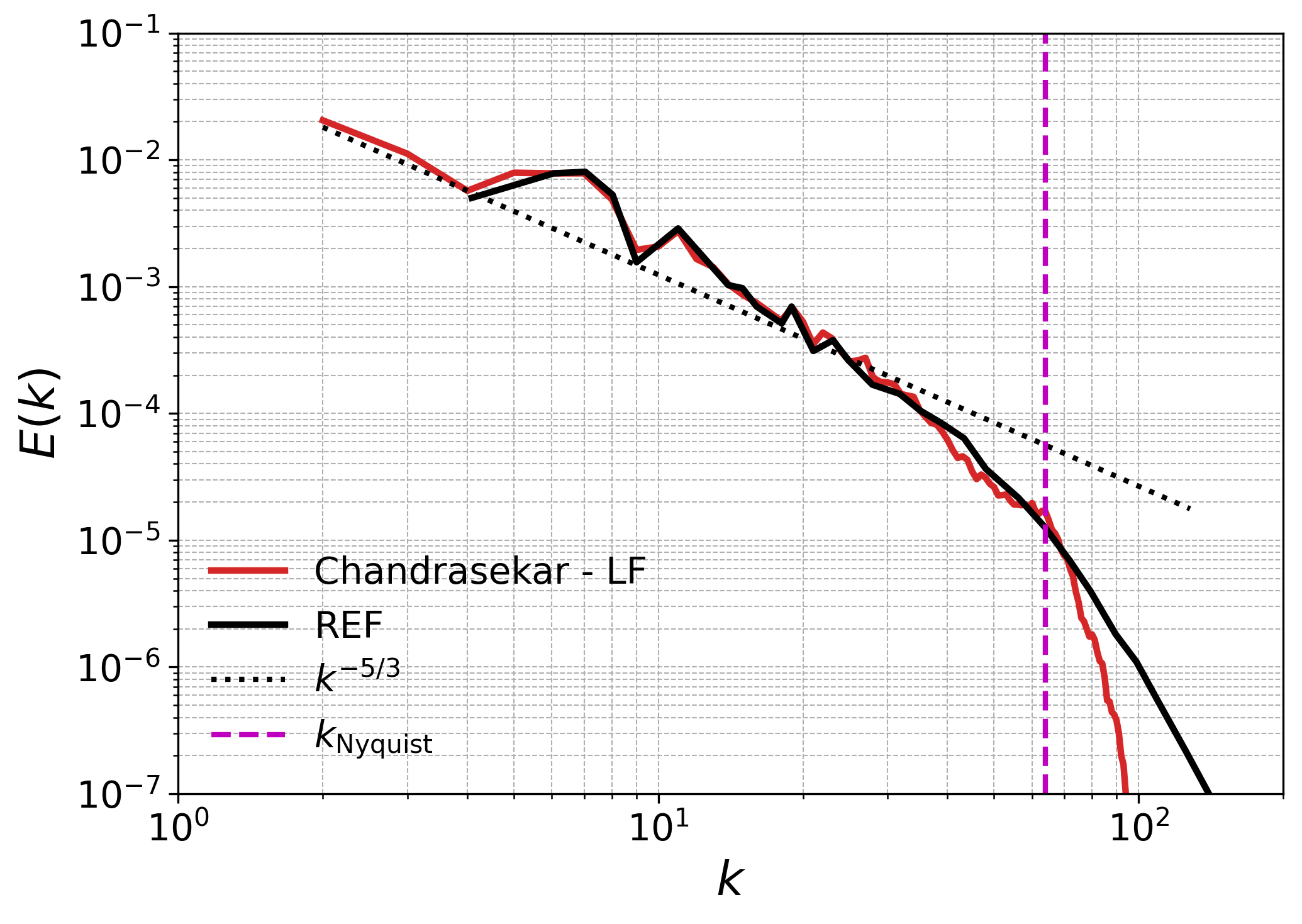}
        \caption{Kinetic energy spectrum at $t/t_c=9$.}
        \label{fig:TGV1600_spectrum}
    \end{subfigure}
    \caption{Validation of the GPU-accelerated \texttt{HORSES3D} solver for the Taylor--Green vortex at $\mathrm{Re}=1600$. 
    Panel (a) shows the kinetic energy dissipation rate compared with the high-resolution DRP reference solution of Bull and Jameson~\cite{Bull2015TGV}. 
    Panel (b) shows the kinetic energy spectrum at $t/t_c=9$ compared with the pseudo-spectral reference solution of Carton de Wiart \textit{et al.}~\cite{CartonDeWiart2014DG}.}
    \label{fig:TGV1600_validation}
\end{figure}

\paragraph{Results:}
Figure~\ref{fig:TGV1600_dissipation} shows the time evolution of the kinetic energy dissipation rate and compares the present results to the high-resolution dispersion-relation-preserving (DRP) reference solution on a fine $512^3$ grid reported by Bull and Jameson~\cite{Bull2015TGV}. The GPU-accelerated \texttt{HORSES3D} solution accurately captures the characteristic dissipation peak associated with the laminar--turbulent transition, as well as the subsequent decay in the fully turbulent regime. Minor discrepancies during the peak-dissipation phase are consistent with the reduced spatial resolution and are typical of high-order implicit LES approaches.
The kinetic energy spectrum at $t/t_c=9$ is shown in Fig.~\ref{fig:TGV1600_spectrum} and compared with the pseudo-spectral reference solution on a $512^3$ grid from Carton de Wiart \textit{et al.}~\cite{CartonDeWiart2014DG}. The spectrum is in good agreement with the reference solution over the energy-containing and intermediate wavenumber ranges, with only minor discrepancies near the cut-off wavenumber, as expected for the present resolution, indicating that the selected numerical configuration provides a robust representation of the flow dynamics.

In general, the results demonstrate that the GPU-accelerated implementation of \texttt{HORSES3D} reproduces the established reference behaviour of the Taylor--Green vortex at $\mathrm{Re}=1600$. 

The agreement in both kinetic energy dissipation and spectral content confirms the correctness of the numerical formulation and validates the GPU port prior to the large-scale performance and scalability studies presented in the following sections.

\subsection{Validation: Turbulent channel flow with wall modelling}
\label{sec:validation_channel}

The GPU implementation of a wall-function boundary condition in \texttt{HORSES3D} is validated using a high-Reynolds-number turbulent channel flow. This test assesses both the correctness of the wall modelling and its robustness under distributed-memory parallel execution.
The wall model follows the Reichardt law-of-the wall formulation following\cite{reichardt1951vollstandige,spalding1961single,frere2017application}, see the implementation details in the Appendix \ref{app:wall_model}. The wall shear stress is reconstructed from a sampled off-wall velocity and is imposed weakly through the viscous numerical fluxes.

The test case consists of a fully developed turbulent channel flow between two parallel walls separated by a distance $2h$. Periodic boundary conditions are imposed in the streamwise and spanwise directions, with domain lengths $2\pi h$ and $\pi h$, respectively. The friction Reynolds number is fixed to $Re_{\tau}=5200$ through a constant pressure-gradient forcing. The Mach number is maintained at $M \approx 0.1$ to approximate incompressible conditions. The DNS database of Lee and Moser~\cite{Lee_2015} is used as reference.

\paragraph{Numerical setup:}
The computational mesh comprises $20\times10\times10$ uniformly distributed hexahedral elements of polynomial order $P=4$. The compressible Navier--Stokes equations are solved using a DGSEM discretisation with Gauss--Lobatto nodes. The inviscid terms are discretised using the Pirozzoli split form, with interface fluxes computed using the Roe Riemann solver. Viscous terms are discretised using the BR1 scheme. Time integration is performed using an explicit third-order Runge--Kutta method. 

\begin{figure}[h]
    \centering

\centering
\includegraphics[width=0.9\textwidth]{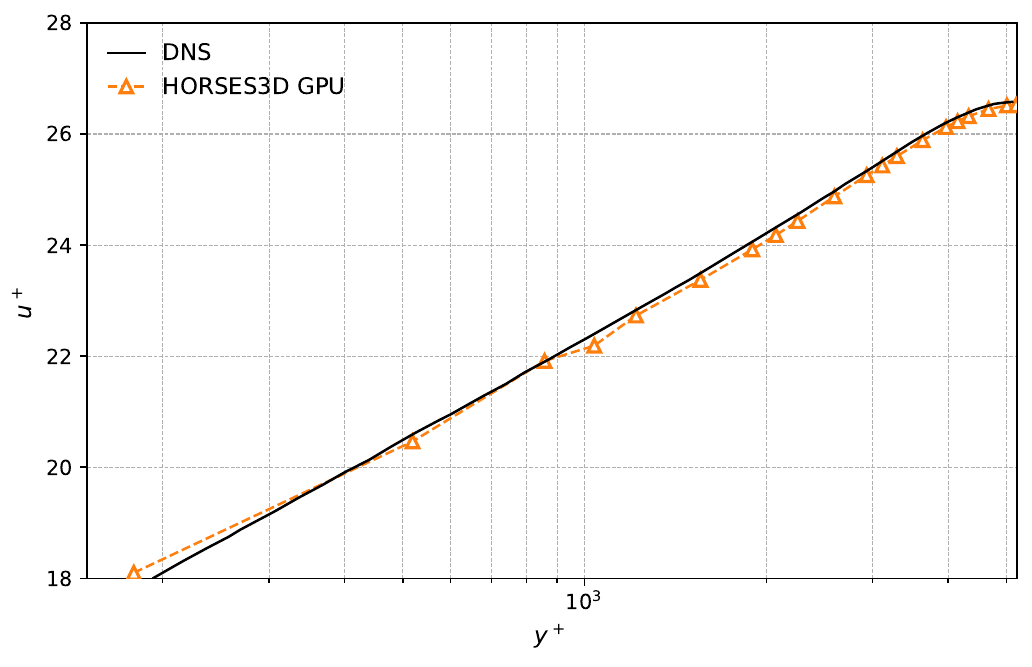}
    \caption{Mean velocity profiles for the turbulent channel flow test case using wall function. Results are compared with reference DNS data (\(Re_\tau = 5200\))~\cite{Lee_2015}.} 
    \label{fig:comparison}
\end{figure}

\paragraph{Results:}
Figure~\ref{fig:comparison} shows the mean velocity profiles in wall units obtained with the GPU implementation. Excellent agreement is observed with respect to the DNS data, confirming that the GPU wall-function implementation reproduces the expected behaviour. 


\section{Performance Analysis of \texttt{HORSES3D-GPU}}
\label{sec:Performance}

\subsection{Hardware}
The performance analysis is performed in MareNostrum 5 ACC, thanks to various EuroHPC and Red Española de Supercomputación (RES) projects. The MareNostrum 5 ACC partition contains 1,120 accelerated compute nodes based on Intel Sapphire Rapids processors and NVIDIA Hopper GPUs. Each node has two CPU sockets, with one Intel Xeon Platinum 8460Y+ processor per socket. Each processor provides 40 cores at 2.3 GHz, resulting in a total of 80 CPU cores per node. In addition, each node contains four NVIDIA Hopper H100 GPUs, each equipped with 64 GB of HBM2e memory. Note that these GPUs use HBM2e rather than the higher-bandwidth HBM3 memory available in other H100 configurations. 
The latter is important, as profiling of the GPU version of \texttt{HORSES3D} indicates that multiple performance critical kernels are memory bound. Therefore, the available bandwidth has a direct effect on the code performance. Each node has sixteen DIMM 32GB 4800MHz DDR5 for a total of 512GB of RAM. Intranode communications are through NVLink 2.0 and internode communications are through four ConnectX-7 NDR200 InfiniBand for a total of 800Gbit/s of bandwidth per node. MareNostrum 5 is one of the most advanced supercomputing facilities in the world with a theoretical peak performance of 260 PFlops for its accelerated partition.

\subsection{Performance analysis - Taylor--Green Vortex}

Our first benchmark uses the TGV problem described in Section \ref{sec:validation_TGV}. The benchmarking has been conducted using a Taylor Green Vortex (TGV) test case using a cubic domain split into non--overlapping hexahedral elements. The smallest mesh in the numerical experiments comprises 4,096 elements ($16\times16 \times 16$) and the largest of 4,194,304 elements ($256 \times 128 \times 128$). We used polynomial orders from P=3 to P=7 to estimate the performance and efficiency of the code for a wide spectrum and find which configuration makes the best use of available resources. 

\texttt{HORSES3D} offers a plethora of numerical setups to perform LES simulations in GPUs, as detailed in Section~\ref{sec:gpu_port}. 
Here,
we compare on a standard DGSEM (SDG) configuration using Gauss-Legendre (LG) points, without explicit subgrid model, and on a kinetic energy preserving scheme (KEP) that uses Kennedy and Grubber fluxes on Gauss-Lobatto points. The latter has been found to be robust and accurate for complex test cases, often without the need for an explicit subgrid model \cite{NTOUKAS2025104425,rubio2025explicitsubgridmodelsenhance} . 

The performance of the code is quantified using the performance index (PID) metric, a commonly used metric for measuring computational efficiency \cite{krais2021flexi,Kurz2025}, defined as
\begin{equation}
\text{PID} = \frac{\text{Walltime} \times \#\text{Ranks}}{\#\text{RK stages} \times \text{DOF}\times \#\Delta t}.   
\label{eq:eff}
\end{equation}
This metric measures the time required to evaluate the right-hand side (RHS) and advance the solution in time per degree of freedom (DOF). For example, for an explicit $3^{rd}$ order Runge-Kutta scheme, we have three RHS evaluations per non-dimensional time step. 
The $\#\text{Ranks}$ in this work is the number of GPUs used to run the simulation. 
The total number of degrees of freedom is computed as $\text{DOF} = N_e \times (P+1)^3$, where $N_e$ is the number of hexahedral elements. All performance metrics reported here were collected over 1000 non-dimensional time steps, $\#\Delta t$. 

To establish the baseline performance, the code is first executed on a single GPU. The problem size and number of devices are then progressively increased to evaluate both strong and weak scalability. In addition to demonstrating the computational performance and parallel efficiency of the code, these studies provide practical guidance in selecting the problem size and hardware configuration required to achieve optimal performance.

\subsubsection{Single GPU performance}

\begin{figure}[H]
    \centering
    \includegraphics[width=0.9\textwidth]{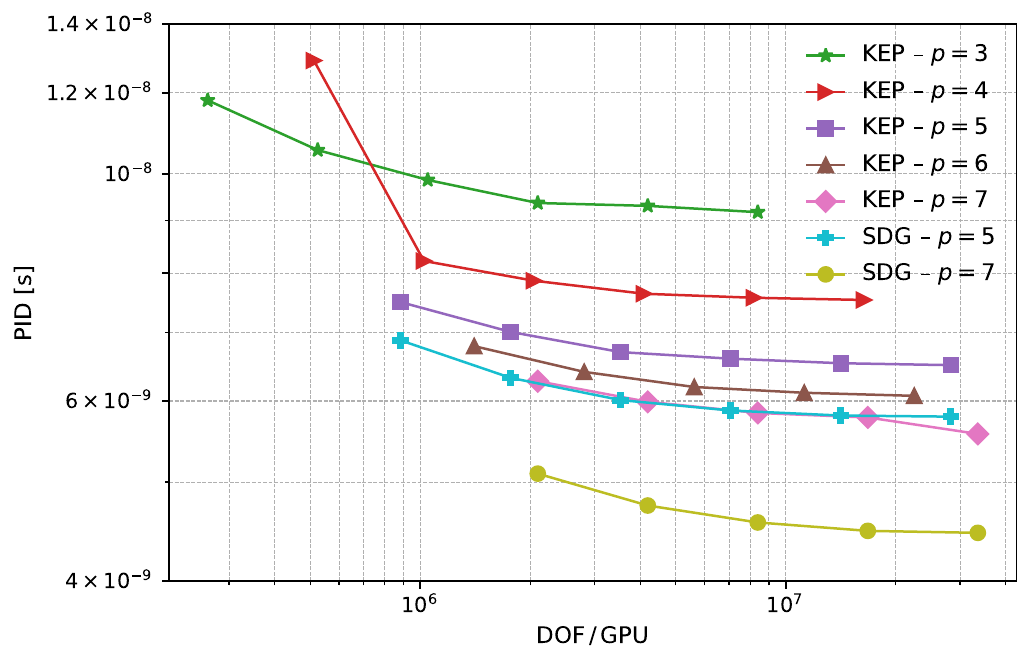}
    \caption{Single GPU performance of \texttt{HORSES3D} for the KEP and SDG schemes for different polynomial orders and different mesh sizes. The figure depicts performance index, PID, (Eq. \ref{eq:eff}, the lower the better) vs degrees of freedom (DOF) per GPU.} 
   \label{fig:singleGPUPerf}
\end{figure}

The first set of tests is performed on a single Nvidia H100 GPU. We use meshes with 4,096 up to 131,072 elements for $P = 3$ to $P = 5$ and 4,096 up to 65,536 elements for $P=6$ and above. This is because for $P>5$, we cannot fit a mesh with 131,072 elements on a single GPU. The results are presented in Figure~\ref{fig:singleGPUPerf} for the standard DGSEM (SDG) and the kinetic energy preserving DGSEM (KEP). The first thing to notice from Figure~\ref{fig:singleGPUPerf} is that the efficiency of the solver increases as we increase the polynomial order, measured via the PID metric. Moving from $P = 3$ to $P = 7$ increases efficiency by a factor of 2. This result is expected given the particular implementation of \texttt{HORSES3D}, the locality of the data, and the structure of the DGSEM operations. Increasing the polynomial order results in more degrees of freedom being processed locally within each element, thereby increasing the amount of work performed by each launched kernel and reducing the relative pressure on device memory bandwidth. This observation is particularly relevant because it indicates that for a fixed number of degrees of freedom, $p$--refinement can be more efficient on GPUs than $h$--refinement. This computational advantage is consistent with a key property of high-order DGSEM discretisations: for sufficiently smooth solutions, increasing the polynomial order can yield exponential convergence, whereas mesh refinement alone generally provides polynomial convergence in $h$. High-order DGSEM discretisations in \texttt{HORSES3D} therefore offer an additional advantage on GPU architectures: improved computational efficiency.

The second important result of Figure~\ref{fig:singleGPUPerf} is that if we increase the load per GPU or the degrees of freedom per GPU (DOF/GPU), we see an increase in the solver efficiency. This is also expected because it is related to the occupancy of computing resources. By offloading more work to the GPU, we can hide the latency of the device. The results indicate that for a load of more than 16384 elements per GPU, we are in $90\%$ of the optimum performance achieved with the highest number of DOFs for each polynomial order. This trend is consistent across all curves and it is a strong indicator for setting up the \texttt{HORSES3D-GPU} solver. More specifically, to make efficient use of the hardware, a minimum of 16,000 to 20,000 elements per GPU are required. The third interesting result of the results presented in Figure~\ref{fig:singleGPUPerf} is that the SDG is more efficient than the KEP setup in \texttt{HORSES3D}. More specifically, the $P=7$ KEP solver is as efficient as the $P=5$ SDG solver and $25\%$ less efficient than the $P=7$ SDG solver. 
This first study provides valuable information about the performance characteristics of \texttt{HORSES3D} on GPUs and establishes an initial guideline for creating optimal setups that make efficient use of computational resources. 

\subsubsection{Multi GPU performance}

To evaluate the performance of \texttt{HORSES3D} in modern HPC systems, we perform numerical experiments to determine both strong and weak scalability characteristics of the code. We use the same TGV case with a cubical domain and run the simulations for 1000 non-dimensional time steps using an explicit RK3 scheme. We focus on the KEP scheme of \texttt{HORSES3D} as this is the preferred configuration for large and complex industrial cases \cite{NTOUKAS2025104425,rubio2025explicitsubgridmodelsenhance}. 

\begin{figure}[h]
    \centering
\includegraphics[width=0.9\textwidth]{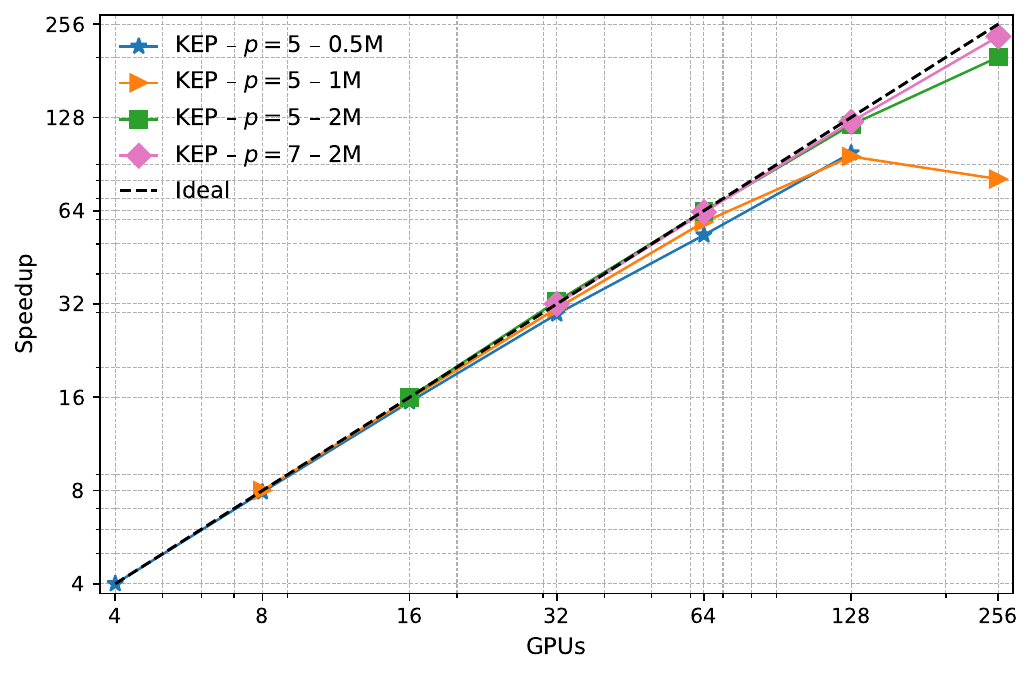}
    \caption{Strong scalability for the KEP scheme for polynomial orders $P=5$ and $P=7$ and meshes with 0.5, 1 and 2 million elements.}
   \label{fig:StrongScalingEC}
\end{figure}

The first set of tests is conducted to determine the strong scalability. The strong scalability is a test measuring the speed-up of the code by keeping the problem size constant and increasing the number of GPUs for each test. Ideally, a code should showcase a speed-up equivalent to the number of GPUs used. The results for the KEP setup are presented in Figure~\ref{fig:StrongScalingEC}. We use three different meshes with 524,288, 1,048,576 and 2097152 elements. We also use $P=5$ and $P=7$, as these are two typical options. The baseline for each curve is set as the performance using the minimum amount of resources that can fit each test case. We then double the number of GPUs and start collecting the performance figures.

The results presented in Figure~\ref{fig:StrongScalingEC} indicate a near-perfect scaling of \texttt{HORSES3D} under the condition that we have a sufficient load per GPU. The data points start to deviate from the ideal scaling curve when the load per GPU falls below 16,000 elements. This is consistent for both $P=5$ and $P=7$, although it is less pronounced for $P=7$. These results indicate that if the load per GPU is less than this threshold, there is not enough work on each GPU to hide the MPI communications. Therefore, some performance loss is expected that can be severe if the load becomes extremely low, as showcased for the case KEP-P5-1M in Figure~\ref{fig:StrongScalingEC}. These results set a guideline for making efficient usage of GPU resources, indicating again that a minimum of 16,000 elements per GPU should be used. This threshold coincides with the threshold for optimum performance on a single GPU device. We observe in Figure~\ref{fig:StrongScalingEC} that if there is sufficient load per GPU, \texttt{HORSES3D} is able to efficiently exploit hundreds of GPUs.

\begin{figure}[H]
    \centering
    \begin{subfigure}{0.48\textwidth}
\includegraphics[width=\textwidth]{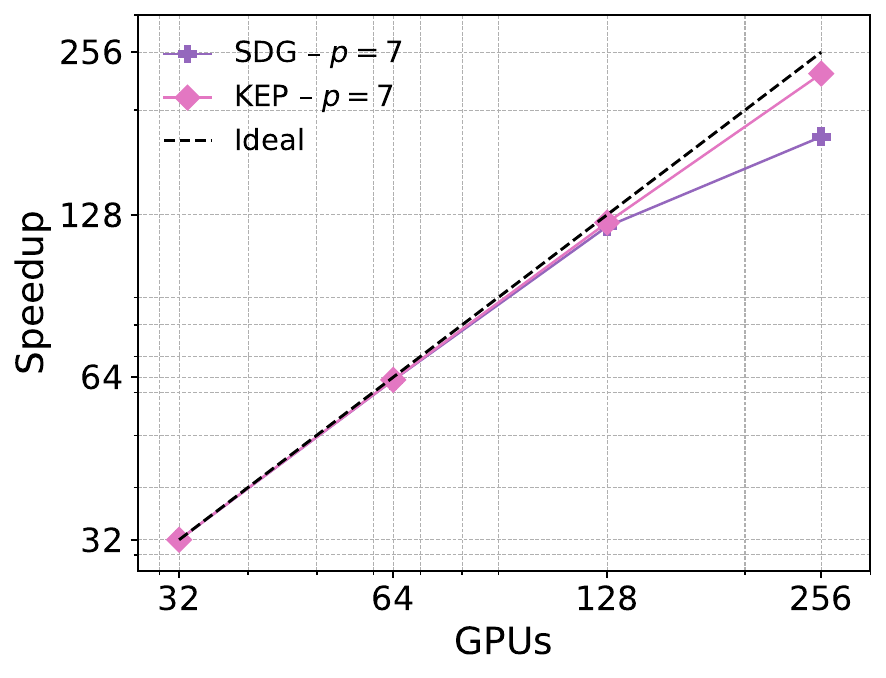}
        \caption{Strong scalability for $P=7$ and 2 million elements.}
        \label{fig:StrongScalingP7}
    \end{subfigure}
    \hfill
    \begin{subfigure}{0.48\textwidth}
        \includegraphics[width=\textwidth]{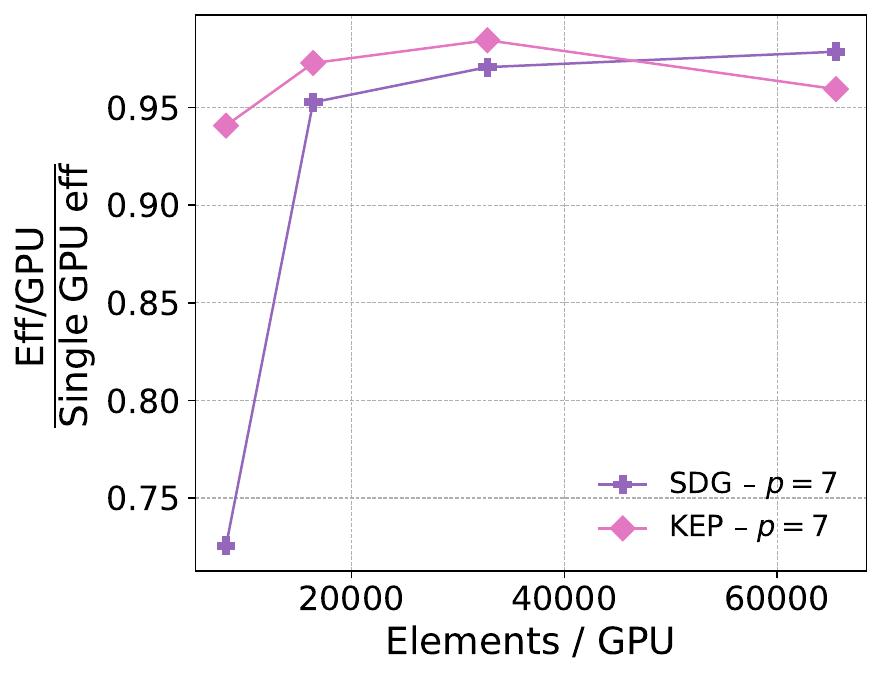}
        \caption{Ratio of Multi-GPU to Single-GPU performance (PID) for specific load per device.}
        \label{fig:MutliToSinglePerf}
    \end{subfigure}
    \caption{Strong scalability results for the KEP and SDG schemes for $P=7$ and a 2 million element mesh.}
    \label{fig:StrongScaling}
    
\end{figure}

Now, we consider the largest case with 2,097,152 elements and $P = 7$, which is approximately 1.07 billion degrees of freedom, and we compare the strong scalability for the KEP and SDG implementations. The smallest number of GPUs that can fit this case is 32 and we test the scaling up to 256 GPUs. The results in Figure~\ref{fig:StrongScalingP7} show that for both schemes we have perfect scaling, and only for a load of 8,192 elements per GPU we start to see performance degradation for the SDG scheme. Also, for the results in Figure~\ref{fig:StrongScalingP7}, the baseline performance is that of 32 GPUs. We also define the ratio of multi to-one GPU simulations, and compare SDG to KEP in Figure~\ref{fig:MutliToSinglePerf}.  For the KEP scheme, we see that for a wide range of elements per GPU, even for lower than 16000 elements, we are close to or above $95\%$ of the performance of a single GPU, showing that \texttt{HORSES3D} is able to hide the MPI communications and scale perfectly for a wide range of loads per GPU. 



\begin{figure}[h]
    \centering
    \includegraphics[width=0.9\textwidth]{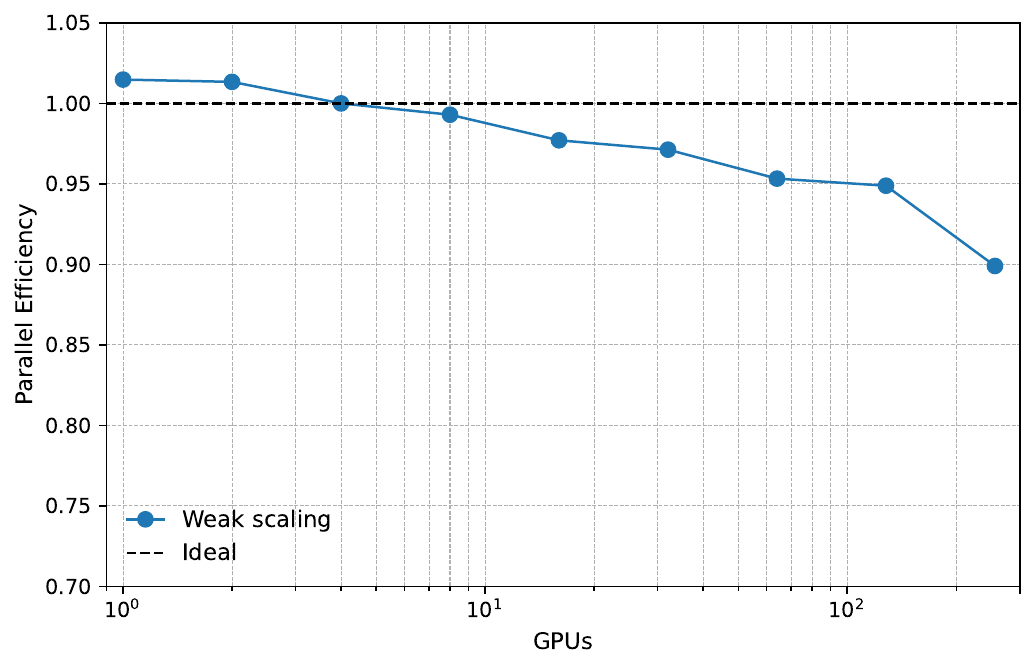}
    \caption{Weak scalability of \texttt{HORSES3D} for $P=5$ and 16384 elements per GPU. Results are relative to single node performance (4 GPUs).}
   \label{fig:Weakscaling}
\end{figure}

The last metric to evaluate the parallel performance of the GPU version of \texttt{HORSES3D} is to evaluate the weak scalability. In this test, we maintain the load per device constant, and we measure the time to complete 1,000 non-dimensional time steps as we increase the number of GPUs. The baseline performance is that for a single node and 4 GPUs. For this series of tests, we select the KEP scheme with $P=5$. The load per GPU is set to 16,384 elements. Thus, the first mesh on a single GPU has 16,384 elements and 3.5 million degrees of freedom, and the last one on 256 GPUs has 4.2 million elements and 905 million degrees of freedom. The results are presented in Figure~\ref{fig:Weakscaling}. For a perfect weak scaling the points should be above or as close as possible to 1. We observe that for 1 and 2 GPUs we obtain slightly better performance, which is expected. As we increase the number of nodes and GPUs we see that the performance drop slightly. This is because the work is spread out on GPUs (which are not physically adjacent in Marenostrum) and internode MPI communications become less efficient, therefore reducing the performance of \texttt{HORSES3D}. This particularity cannot be controlled by the user in most cases and within heavily used HPC facilities. Nonetheless, the code runs efficiently under these conditions. The results indicate that for 128 GPUs or less we are in the $95\%$ of the single-node performance and in the $90\%$ for 256 GPUs. These results show that \texttt{HORSES3D} can run efficiently and takes advantage of state-of-the-art hardware. These figures were obtained for $P=5$ and for the minimum load per GPU, as indicated in the results in this work. We expect that for higher number of elements per GPU and higher polynomial order, the weak scalability to be even better for larger number of nodes and GPUs.

\section{Towards aeronautic simulations on GPUs - Common Research Model}
\label{sec:CRM}

To demonstrate the capabilities of the GPU-accelerated \texttt{HORSES3D} solver in industrially relevant geometries, we present strong and weak scalability results for the High-Lift Common Research Model (CRM-HL)~\cite{clark2025hlpw}. The CRM-HL is an open, publicly available geometry representative of a commercial transport aircraft in high-lift configuration and has been widely used for CFD validation, notably in the AIAA High-Lift Prediction Workshops, including HLPW-5. This section is carried out with the test Case 1 which contains only the wing and body geometry. 

\subsection{Computational Setup}

The simulations are performed using \texttt{HORSES3D} and the compressible Navier–Stokes formulation, employing Chandrashekar entropy-conserving averaging together with Gauss–Lobatto points. The Lax–Friedrichs Riemann solver is used for numerical fluxes, while the viscous terms are discretised using BR1 formulation. This combination of entropy-conserving averaging and flux discretisation has been found to be  robust for challenging simulations.
Time is discretised with an RK3 scheme. 
The Mach number is $M=0.2$ and the Reynolds number of $\mathrm{Re}=5.6\times10^6$,  based on the mean aerodynamic chord of the wing and an angle of attack of $\mathrm{AoA}=11^\circ$. All simulations are run for 1,000 non-dimensional time steps to measure the efficiency of the PID, defined in Eq.~\eqref{eq:eff}.
All timing measurements exclude mesh generation and storage operations; only the residual evaluation phase is considered, as recorded by the internal stopwatch instrumentation.

Four levels of hexahedral mesh refinement are used, labelled h0 through h3, with element counts ranging from approximately 3.3 million (h0) to 20.8 million (h3) elements, as detailed in Table~\ref{tab:CRM_grids}. 
\begin{table}[h!]
\centering
\small
\caption{Number of hexahedral elements for each grid.}
\begin{tabular}{lrr}
\hline
\textbf{Grid} & \textbf{Elements}  \\
\hline
h0 &  3,342,208  \\
h1 &  5,839,712  \\
h2 & 11,658,356  \\
h3 & 20,833,704  \\
\hline
\end{tabular}
\label{tab:CRM_grids}
\end{table}

\begin{figure}
    \centering
    \includegraphics[width=0.6\linewidth]{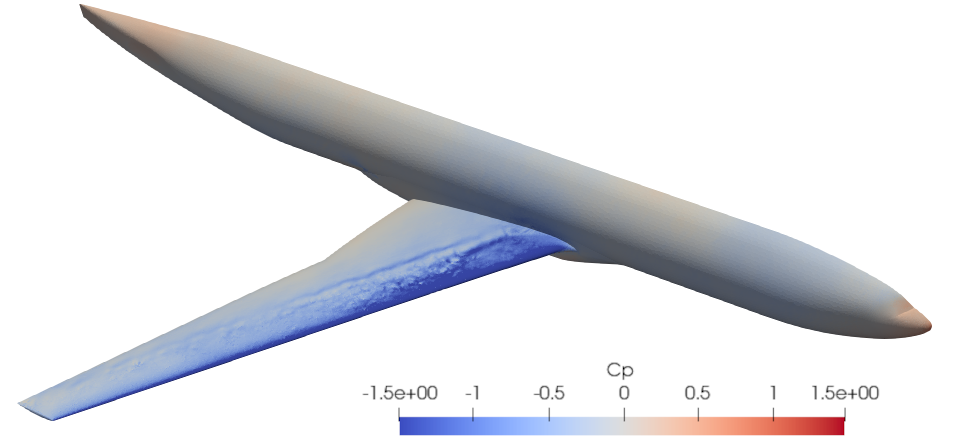}
    \caption{Representative flow solution used in the performance assessment. Instantaneous pressure coefficient distribution on the CRM configuration obtained on the h1 grid with $P=2$ ($\approx 157.7\times10^6$ degrees of freedom).}
    \label{fig:CRM}
\end{figure}

The range covered by the benchmark is substantial: the smallest configuration (h0, $P=1$) contains approximately $26.7$ million DOF, while the largest (h3, $P=7$) reaches approximately $10.7$ billion DOF, spanning nearly three orders of magnitude and enabling a thorough assessment of the solver across a wide range of problem sizes.

Figure~\ref{fig:CRM} presents a representative solution obtained on the h1 grid using a polynomial order of $P=2$, corresponding to approximately $157.7$ million degrees of freedom. 
The pressure coefficient distribution illustrates the complex aerodynamic loading over the CRM configuration and demonstrates the capability of the solver to handle realistic industrial-scale flow simulations.

\subsection{Performance analysis}

The scaling studies have been carried out on the MareNostrum 5 accelerated partition, using up to 2048 GPUs, corresponding to half of the full system capacity (4,096 GPUs). 
%

\subsubsection{Baseline performance and scaling strategy}
 Due to the significantly larger problem sizes considered in the present configurations, single-GPU executions are not feasible in practice. For this reason, all performance measurements are obtained using multi-GPU runs.

To enable a consistent comparison across different configurations, the results are normalised with respect to the number of GPUs and reported in terms of the efficiency metric PID defined in Eq.~\eqref{eq:eff}. While this rescaling allows a coherent comparison between setups, it must be emphasised that perfect linear scaling is not achieved in practice due to communication overheads and load imbalance effects. Nevertheless, the resulting metrics provide a meaningful basis for assessing the relative efficiency of the different discretisations and mesh resolutions.

\begin{figure}[h!]
\centering
\includegraphics[width=0.9\textwidth]{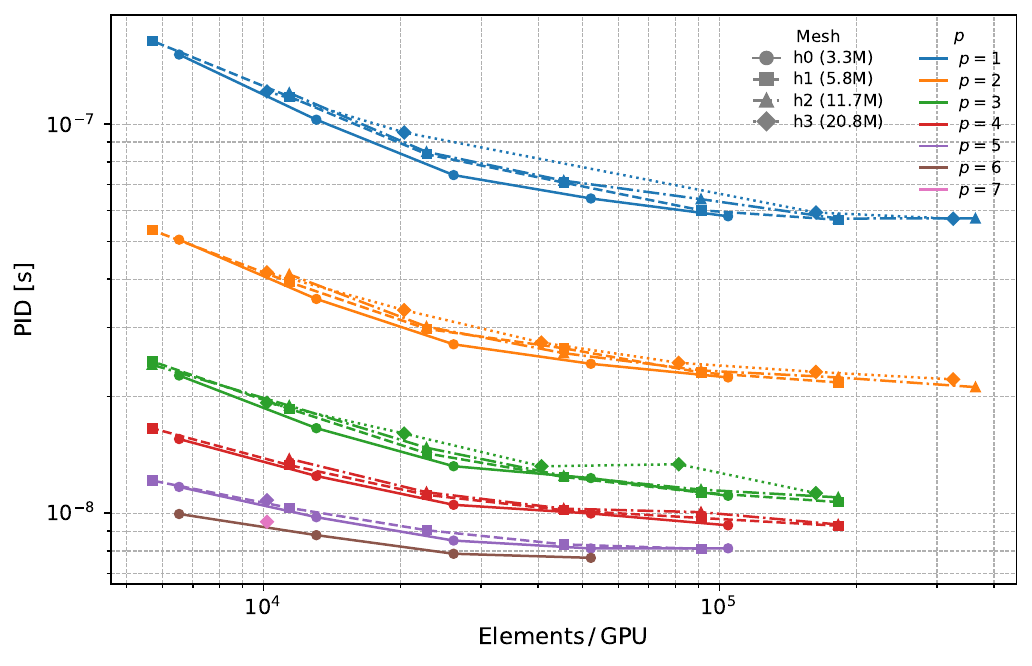}
\caption{Performance index PID (Eq.~\eqref{eq:eff}) for the CRM wing–body configuration across mesh levels h0--h3 and polynomial orders $P=1$--$7$. The results are shown as a function of elements per GPU, highlighting the onset of efficiency saturation once a minimum load of approximately $2\times10^4$ elements per device is reached.}
\label{fig:CRM_eff_dof_elem}
\end{figure}

The performance metric PID is first introduced as in the TGV case. The resulting performance index PID is summarised in Figure~\ref{fig:CRM_eff_dof_elem}, which shows the efficiency between mesh levels h0--h3 and polynomial orders $P=1$--$7$ as a function of the number of elements per GPU. As observed, general trends are consistent with those reported for the TGV benchmark, with only moderate variations due to the increased complexity of the CRM configuration. In particular, once a minimum load of approximately $20,000$ elements per GPU is reached, the PID metric becomes essentially independent of the problem size, indicating that the computation operates in a regime of near-optimal GPU utilisation.

\begin{figure}[h!]
\centering
\includegraphics[width=0.9\textwidth]{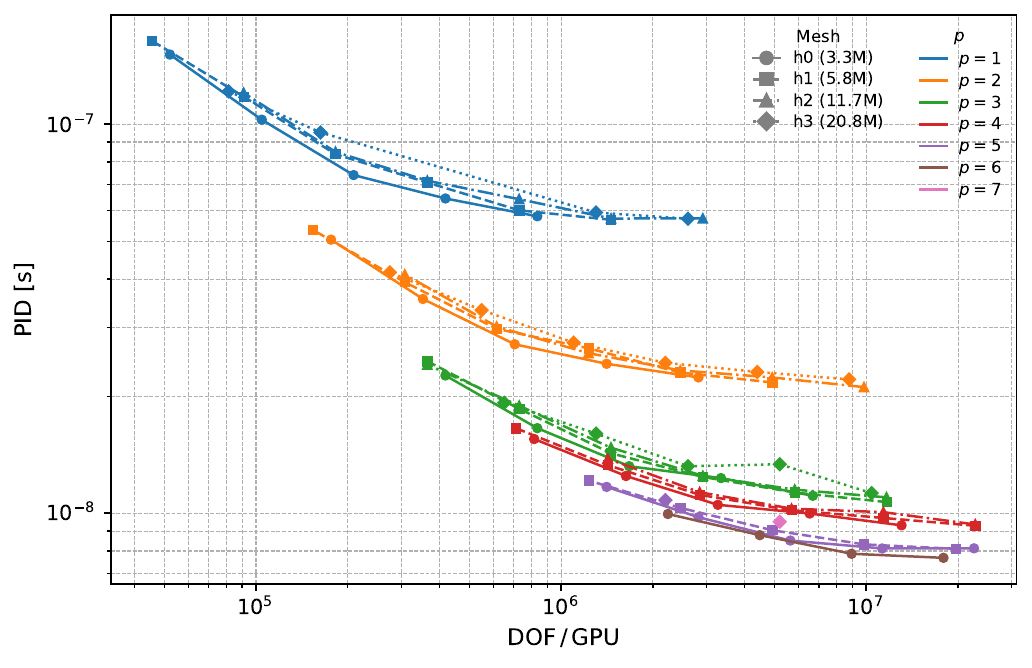}
\caption{Performance index PID (Eq.~\eqref{eq:eff}) for the CRM wing–body configuration across mesh levels h0--h3 and polynomial orders $P=1$--$7$. The results are reported as a function of degrees of freedom per GPU, allowing direct assessment of efficiency saturation with increasing load per device.}
\label{fig:CRM_eff_dof}
\end{figure}

For completeness, Figure~\ref{fig:CRM_eff_dof} provides the same efficiency metric reorganised by mesh level and polynomial order as a function of the number of DOF per GPU, enabling a more direct comparison with the single-GPU TGV baseline case. This representation confirms that the trends with respect to polynomial order are preserved in the CRM configuration, with higher orders consistently yielding improved efficiency due to increased computational intensity per degree of freedom.
Direct comparison with the single-GPU TGV results reported in Section~\ref{sec:Performance} confirms that the CRM efficiency values remain of the same order of magnitude. In particular, for the highest polynomial orders ($P=5$--$7$), the CRM configuration achieves $\mathrm{PID} \approx 8$--$10 \times 10^{-9}$~s at optimal load, which closely matches the single-GPU TGV performance of $\mathrm{PID} \approx 4$--$6 \times 10^{-9}$~s. For lower orders ($P=3$), the CRM efficiency is approximately $\mathrm{PID} \approx 1 \times 10^{-8}$~s, again consistent with the TGV baseline of $\mathrm{PID} \approx 9 \times 10^{-9}$~s.

The slightly higher values observed for the CRM configuration are attributed to two main factors. First, the CRM geometry introduces additional computational overhead associated with complex boundary condition treatment (including viscous walls, inlet, and outlet conditions) that is absent in the periodic TGV setup. Second, the reported efficiency is obtained from wall-clock measurements of multi-GPU runs and rescaled under the assumption of ideal strong scalability; consequently, any loss of parallel efficiency due to MPI communication, load imbalance arising from unstructured mesh partitioning, or network contention manifests itself as an apparent increase in PID relative to the single-GPU baseline.
Despite these effects, the close agreement between CRM and TGV efficiency levels demonstrates that the GPU-accelerated \texttt{HORSES3D} solver preserves its performance characteristics when applied to geometrically complex, industrially relevant configurations.

\subsubsection{Strong scalability results}

\begin{figure}[h!]
    \centering
    \includegraphics[width=1\textwidth]{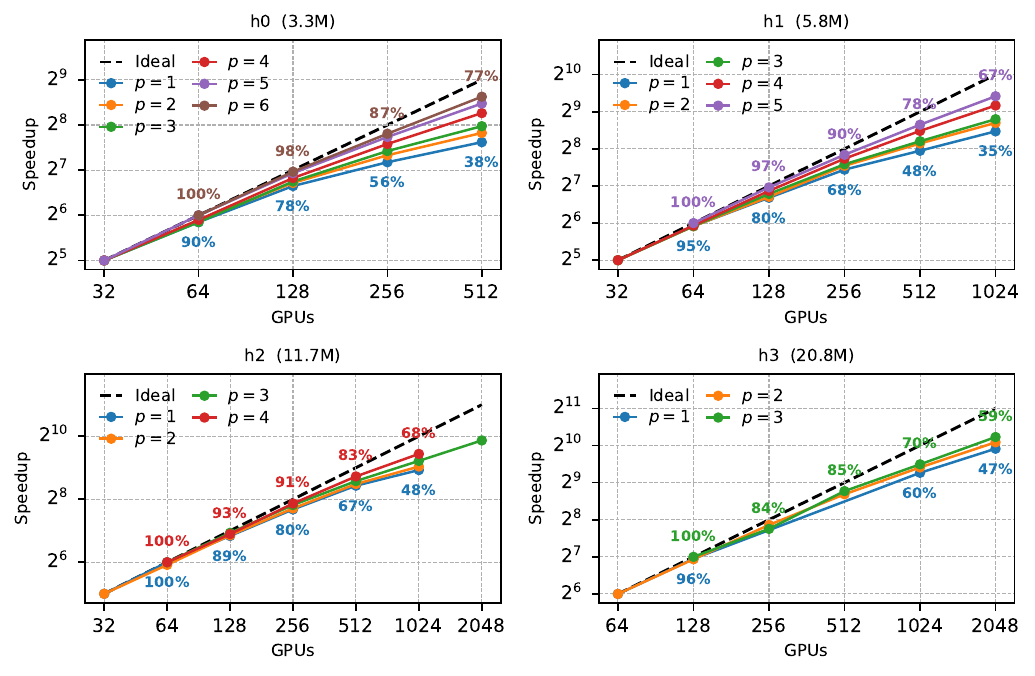}
    \caption{Strong scaling speedup relative to the baseline GPU count for the CRM wing–body configuration. Results are shown for mesh levels h0 (3.3M hexas), h1 (5.8M hexas), h2 (11.7M hexas), and h3 (20.8M hexas).}
    \label{fig:CRM_speedup}
\end{figure}
Regarding the strong scaling behaviour shown in Figure~\ref{fig:CRM_speedup}, the solver exhibits the expected dependence on the load per GPU across all mesh levels. Overall, near-optimal efficiency is achieved when the number of elements per GPU remains above the previously identified threshold of $\mathcal{O}(10^4)$.
For the h1 mesh ($\approx 5.8$ million elements), strong scaling remains close to ideal for moderate GPU counts. For example, a speedup of $1.91$ is obtained when doubling the number of GPUs from 32 to 64 for $P=4$, with approximately $9.1\times10^4$ elements per GPU. As the number of GPU increases to 1024, the load decreases to approximately $5.7\times10^3$ elements per GPU, and the efficiency decreases to $56\%$, indicating the onset of communication-dominated behaviour.
A similar trend is observed for $P=5$ in the same mesh, although with improved robustness at intermediate loads. At 256 GPUs, corresponding to approximately $2.3\times10^4$ elements per GPU, the efficiency remains high ($\sim 89.6\%$), while at 1024 GPUs it reduces to $66.8\%$ as the workload per-device falls below the optimal regime. This confirms the increased resilience of higher polynomial orders to reduced per-GPU load, consistent with the TGV observations.
The largest configuration (h3, $\approx 20.8$ million elements) is used to assess scalability at extreme problem sizes, reaching up to 2048 GPUs. In this regime, the solver maintains strong scaling behaviour close to optimal up to approximately 1024 GPUs, where the load is still around $2.0\times10^4$ elements per GPU and efficiencies of the order of $80\%$ are observed for $P=2$.
When further increasing the number of GPUs to 2048, the workload decreases to approximately $1.0\times10^4$ elements per GPU, and the efficiency drops to around $53\%$. This transition is consistent with entering a regime where inter-node communication and global synchronisation start to dominate the computational cost.
It should be noted that the largest execution reported in this study corresponds to the h3 mesh with $P=7$ on 2048 GPUs, resulting in approximately $10.7\times10^9$ degrees of freedom. 
Given the recent commissioning of the system (2023), this places the present study among the early CFD applications exploiting a significant fraction of its GPU resources.
Overall, these results confirm that the GPU-accelerated \texttt{HORSES3D} solver is capable of sustaining efficient execution for geometrically complex, industrial-scale CFD problems at the billion-degree-of-freedom level, while maintaining predictable scalability behaviour across modern HPC architectures.


\subsubsection{Weak scalability results - Load-dependent performance clustering}

To further analyse the dependence of solver performance on the computational load per GPU, a clustering-based weak scaling study is performed. The results are obtained from the entire CRM dataset by grouping configurations with similar elements-per-GPU ratios using a tolerance-based clustering strategy. More specifically, configurations are assigned to the same cluster if their elements-per-GPU value differ by less than $25\%$. This procedure enables the identification of quasi-equivalent load regimes across different mesh levels, polynomial orders, and GPU counts, providing a unified view of performance trends.

\begin{figure}[h!]
    \centering
    \includegraphics[width=0.9\textwidth]{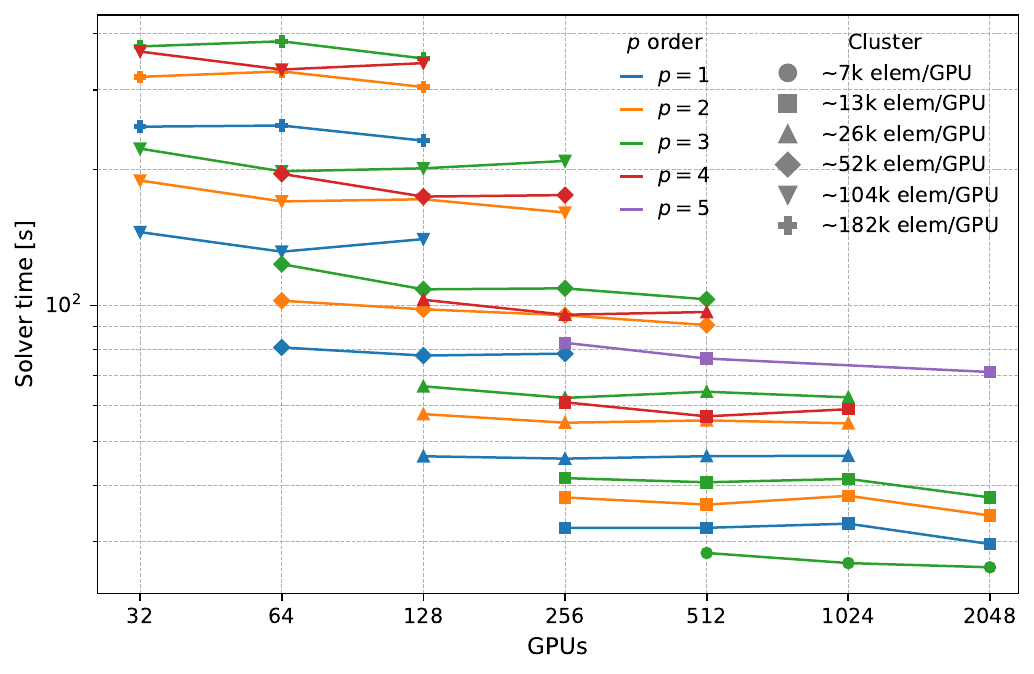}
    \caption{Weak scaling behaviour of the CRM wing–body configuration grouped by clusters of similar elements-per-GPU ratios. Each curve represents a set of configurations with approximately constant computational load per GPU across different mesh levels, polynomial orders, and GPU counts. The results highlight the dominant role of load per device in determining solver efficiency and scaling behaviour.}
    \label{fig:CRM_cluster_scaling}
\end{figure}

The resulting behaviour is shown in Figure~\ref{fig:CRM_cluster_scaling}, where the execution time of the solver is reported as a function of the number of GPUs for each identified cluster. Each curve therefore represents a family of configurations with approximately constant elements-per-GPU load, while varying the global problem size and parallel decomposition.
Across all polynomial orders, the results confirm that performance is primarily governed by the per-GPU workload rather than the absolute number of degrees of freedom. In particular, the curves exhibit near-flat behaviour over a wide range of GPU counts for all configurations, indicating that weak scaling is well preserved across the different load regimes considered.

\section{Conclusions}

In this work, we have presented the GPU-accelerated performance of the high-order DG solver \texttt{HORSES3D} on modern HPC architectures, using both the canonical Taylor–Green vortex (TGV) benchmark and the industrially relevant CRM wing–body configuration. The study provides a comprehensive assessment of solver efficiency, scalability, and robustness across a wide range of polynomial orders, mesh resolutions, and GPU counts.

The results obtained for the TGV benchmark establish the baseline performance characteristics of the solver, showing that efficiency increases with polynomial order due to improved data locality and higher arithmetic intensity. In particular, the analysis demonstrates that performance is governed primarily by the number of elements per GPU, with an optimal regime identified at approximately 16,000 to 20,000 elements per device. This threshold ensures sufficient workload to hide memory latency and communication overheads and is shown to be consistent across all configurations considered.

The CRM configuration extends this analysis to a complex, geometrically realistic case including boundary conditions and non-trivial flow features. Despite the increased complexity, the solver retains performance levels comparable to the TGV benchmark, with efficiency metrics of the same order of magnitude. This shows that the GPU implementation of \texttt{HORSES3D} preserves its computational characteristics when applied to industrial-scale CFD problems.

Strong scaling studies show near-ideal behaviour as long as the load per GPU remains above the identified threshold. As the workload decreases, a gradual degradation in efficiency is observed due to the increasing impact of inter-GPU communication and synchronisation. Weak scaling analysis, including the clustering-based approach, further confirms that the dominant parameter controlling performance is the elements-per-GPU ratio, with consistent behaviour across mesh levels and polynomial orders.

The solver scales efficiently to 2048 GPUs on the MareNostrum 5 system, 
The largest simulation performed in this study reaches approximately $10.7\times10^9$ degrees of freedom, representing the largest computation performed with \texttt{HORSES3D} to date and one of the first CFD applications to take advantage of a substantial fraction of the MareNostrum 5 accelerated partition.

Overall, the results presented in this work confirm that \texttt{HORSES3D} is capable of delivering high accuracy, efficiency and excellent scalability for high-order CFD simulations on state-of-the-art GPU-based HPC systems. The combination of robustness, performance, and flexibility makes it a suitable tool for large-scale industrial applications and positions it as a strong candidate for next-generation exascale computing environments. We have demonstrated that the GPU version of \texttt{HORSES3D} is able to exploit the current generation of HPC architectures and serves as a useful and competitive tool to address the present and future demands of AI.

Future work will investigate the use of higher polynomial orders and larger computational meshes to further exploit modern GPU architectures. Particular emphasis will be placed on production-scale aerodynamic applications, where the favourable performance characteristics observed for high-order discretisations may translate into significant reductions in time-to-solution.



Finally, \texttt{HORSES3D} includes a wide range of multiphysics capabilities that have not been explored in this work, such as multiphase flows, shock-capturing techniques, actuator line models, immersed boundary methods, and aeroacoustic simulations. Parts of these modules have already been ported to GPUs, and their performance and applicability to large-scale simulations will be the subject of future work.

\section*{Acknowledgments}
The authors
acknowledge the funding received by the Grant DeepCFD (Project No. PID2022-137899OB-I00) funded by MICIU/AEI/10.13039/501100011033 and by ERDF, EU and also the funding from the European Union (ERC, Off-coustics, project number 101086075).
This project has received funding from the European Research Council (ERC) under the European Union's Horizon Europe research and innovation programme (grant agreement No. 101167322 - TRANSDIFFUSE).
Views and opinions expressed are, however, those of the authors only and do not necessarily reflect those of the European Union or the European Research Council. Neither the European Union nor the granting authority can be held responsible for them.
The authors gratefully acknowledge the EuroHPC JU Benchmark Access project EHPC-BEN-2024B11-035 and the RES projects RES-IM-2025-2-0013 and RES-IM-2025-3-0007 for providing access to the MareNostrum 5 supercomputer, which enabled the simulations presented in this work. The authors also acknowledge the EPICURE program under the EuroHPC project EHPC-DEV-2025D12-017 for its support in profiling and optimizing \texttt{HORSES3D} on GPU architectures, as well as for organizing the 2026~High-Scalability Workshop on MareNostrum 5.
Andrés M. Rueda-Ramírez gratefully acknowledges funding from the Spanish Ministry of Science, Innovation, and Universities through the ``Beatriz Galindo'' grant (BG23-00062).
Andrés M. Rueda-Ramírez acknowledge funding through the German Federal Ministry for Education and Research (BMBF) project ``ICON-DG'' (01LK2315B) of the ``WarmWorld Smarter'' program.

\appendix
\section{Compressible Navier--Stokes solver}
\label{sec:cNS}
 In this work, we solve the 3D Navier--Stokes equations that can be compactly written as: 
%
\begin{equation}
\boldsymbol{q}_t+ \nabla \cdot {\ssvec{F}}_e = \nabla\cdot\ssvec{F}_{v,turb},
\label{eq:compressibleNScompact}
\end{equation}
where $\boldsymbol{q}$ is the state vector of large scale resolved conservative variables $\boldsymbol{q} = [ \rho , \rho v_1 , \rho v_2 , \rho v_3 , \rho e]^T$, $\ssvec{F}_e$ are the inviscid, or Euler fluxes,
\begin{equation}
\ssvec{F}_e = \left[\begin{array}{ccc} \rho v_1 & \rho v_2 & \rho v_3 \\
                                                                                \rho v_1^2 + p & \rho v_1v_2 & \rho v_1v_3 \\
                                                                                	\rho v_1v_2 & \rho v_2^2 + p & \rho v_2v_3 \\
                                                                                	\rho v_1v_3 & \rho v_2v_3 & \rho v_3^2 + p \\
                                                                                	\rho v_1 H & \rho v_2 H & \rho v_3 H
\end{array}\right],
\end{equation}
where $\rho$, $e$, $H=e+p/\rho$, and $p$ are the large scale density, total energy, total enthalpy and pressure, respectively, and $\vec{v}=[v_1,v_2,v_3]^T$ is the large scale resolved velocity components. Additionally, $\ssvec{F}_{v,turb}$ defines the viscous and turbulent fluxes,
\begin{equation}
\ssvec{F}_{v,turb}= \left[\begin{array}{ccc}0 & 0 & 0\\
 \tau_{xx} & \tau_{xy} & \tau_{xz} \\
 \tau_{yx} & \tau_{yy} & \tau_{yz} \\
 \tau_{zx} & \tau_{zy} & \tau_{zz} \\
 \sum_{j=1}^3 v_j\tau_{1j} + \kappa T_x& \sum_{j=1}^3 v_j\tau_{2j} + \kappa T_y& \sum_{j=1}^3 v_j\tau_{3j} + \kappa T_z
\end{array}\right],
\label{eq:viscousfluxes}
\end{equation}
where $\kappa$ is the thermal conductivity, $T_x, T_y$ and $T_z$ denote the temperature gradients and the stress tensor $\boldsymbol{\tau}$ is defined as $\boldsymbol{\tau} = \mu(\nabla \vec{v} + (\nabla \vec{v})^T) - 2/3\mu \boldsymbol{I}\nabla\cdot\vec{v}$, with $\mu$ the dynamic viscosity 
 and $\boldsymbol{I}$ the three-dimensional identity matrix. 
%


\section{Wall-model formulation and implementation details}
\label{app:wall_model}

\paragraph{Reichardt wall law.}

The wall model is based on the Reichardt law-of-the wall~\cite{frere2017application}, given by
\begin{equation}
u^+_{\parallel} =
\frac{1}{\kappa}\ln\left(1+\kappa y^+\right)
+
\left(
C-\frac{1}{\kappa}\ln\kappa
\right)
\left(
1-e^{-y^+/11}-\frac{y^+}{11}e^{-y^+/3}
\right),
\end{equation}
where $\kappa$ and $C$ are model constants.

The friction velocity is obtained by solving this nonlinear relation locally at each quadrature point.

\paragraph{Velocity projection and wall shear stress.}

The tangential velocity is computed by projecting the sampled velocity onto the wall-parallel direction:
\begin{equation}
\mathbf{x}_{\parallel} =
\frac{
\mathbf{u}_{LES} - (\mathbf{u}_{LES}\cdot\mathbf{n})\mathbf{n}
}{
\left|
\mathbf{u}_{LES} - (\mathbf{u}_{LES}\cdot\mathbf{n})\mathbf{n}
\right|
},
\end{equation}

\begin{equation}
u_{\parallel} = \mathbf{u}_{LES}\cdot\mathbf{x}_{\parallel},
\end{equation}

where $\mathbf{u}_{LES}$ is the sampled velocity and $\mathbf{n}$ is the wall-normal vector. The resulting wall shear stress is projected back to the global reference frame and imposed via the viscous flux formulation.

\paragraph{Local wall-model implementation.}

The implementation uses only locally available data within the wall-adjacent element. The input of the wall-model is extracted from the face opposite to the wall inside the first off-wall element, avoiding dependencies on neighboring elements or MPI partitions.

\paragraph{MPI and Runge--Kutta considerations.}

Due to the ordering of halo exchanges and Runge--Kutta updates, the wall-model input may correspond to solution values from the previous Runge--Kutta stage in some MPI configurations. This avoids additional synchronisation and communication overhead, at the cost of a slight temporal lag.

\paragraph{Additional numerical details.}

The flow is initialized from uniform thermodynamic fields and a laminar parabolic velocity profile. The simulation is advanced until statistically stationary conditions are reached (approximately $25$ flow-through times), followed by averaging over more than $15$ flow-through times. Mean velocity profiles are obtained through temporal and homogeneous spatial averaging.

Additional simulations using different MPI partitionings were performed to assess robustness in distributed-memory environments.